\documentclass[reqno,9 pt]{amsart}
\usepackage{amsmath,amsxtra,amssymb,latexsym, amscd,amsthm}
\usepackage[unicode]{hyperref}
\usepackage{array,tabularx,longtable,multicol,indentfirst,fancybox,color}
\usepackage{graphicx}
\usepackage{multicol}
\usepackage{mathrsfs}
\usepackage{enumerate}
\usepackage{tikz}
\usepackage{tikz-cd} 
\usepackage{tikz,tkz-tab}
\usetikzlibrary{calc, angles, quotes, intersections, positioning, patterns}

\advance\hoffset-1truecm\relax

\newtheorem{theorem}{Theorem}[section]
\newtheorem{definition}{Definition}
\newtheorem{lemma}{Lemma}[section]
\newtheorem{corollary}{Corollary}[section]
\newtheorem{remark}{Remark}
\newtheorem{example}{Example}[section]
\numberwithin{equation}{section}
\def\N{\mathbb{N}}
\def\R{\mathbb{R}}

\def\e{\ell}
\def\p{\partial}
\def\vp{\varphi}
\def\a{\alpha}
\def\b{\beta}
\def\ga{\gamma}

\def\A{\mathscr{A}}

\def\L{\mathcal{L}}
\def\T{\mathscr{T}}
\def\P{\mathcal{P}}
\def\lv{\lVert}

\def\M{\mathcal{M}}

\numberwithin{equation}{section}

\begin{document}
	\title[New polyconvolution product for Fourier-cosine and Laplace integral operators and their applications]{New polyconvolution product for Fourier-cosine and\\ Laplace integral operators and their applications}
	\author[Trinh Tuan]{Trinh Tuan}
	\thanks{Dedicated to the memory of my father (1933–2003).} 
	\thanks{First published: 12 October 2023  \url{https://doi.org/10.1002/mma.9716}}
	\thanks{Received: 5 February 2023 $\big|$ Revised: 15 September 2023 $\big|$ Accepted: 24 September 2023  \textit{Mathematical Methods in the Applied Sciences}}
	
	\thanks{Funding information: This research received no specific grant from any funding agency.}

	\maketitle	
\begin{center}
Department of Mathematics, Faculty of Sciences, Electric Power University,\\ 235-Hoang Quoc Viet Road, Hanoi, Vietnam.\\
E-mail: \textit{tuantrinhpsac@yahoo.com}
\end{center}
\begin{abstract}
The goal of this paper is to introduce the notion of polyconvolution for Fourier-cosine, Laplace integral operators, and its applications. The structure of this polyconvolution operator and associated integral transforms are investigated in detail. The Watson-type theorem is given, to establish necessary and sufficient conditions for this operator to be isometric isomorphism (unitary) on $L_2 (\R_+)$, and to get its inverse represented in the conjugate symmetric form. The correlation between the existence of polyconvolution with some weighted spaces is shown, and Young's type theorem, as well as the norm-inequalities in weighted space, are also obtained.  As applications, we investigate the solvability of a class of Toeplitz plus Hankel type integral equations and linear Barbashin's equations with the help of factorization identities of such polyconvolution. Several examples are provided to illustrate the obtained results to ensure their validity and applicability.

\vskip 0.3cm

\noindent \textbf{MSC CLASSIFICATION}  42A38, 44A10, 44A35, 45E10, 45J05, 26D10. 
\end{abstract}
	
\section{Introduction}
In Fourier analysis and operator theory, the theory of convolutions of integral transforms and convolution type operators have been studied for a long time due to their fundamental role in modeling and solving a wide range of mathematical physics problems. Nowadays, there are two main approaches to the constructing of convolutions for integral transforms. The first one is based on the construction of a generalized shift operator (also called generalized translation operator, or generalized displacement operator). Then the classical translation operator (ordinary translation) in the convolution is replaced by the generalized shift operator, and we get a generalized convolution. The generalized shift operators of the Delsarte\cite{delsarte1939}–Levitan\cite{BMLevitan1962}–Povzner\cite{Povzner1948} type are usually used in these constructions. For instance, this approach is often served to introduce the classical convolution for the Hankel transform. Here we use the idea of the second approach, which rests on the works of Valentin A. Kakichev. His convolution constructing method is based on factorization equality and can be applied to most convolutions with arbitrary integral transform appearing already in \cite{kakichev1967convolution}. In the year 1997, Kakichev generalized this approach and introduced the concept of polyconvolution or "generalized convolution"  \cite{Kakichev1997}, more details as follows.

\begin{definition}\label{defbyKakichev}
	Let $U_i (X_i)$ be the space of linear functions that can differ on the same field and $V (Y)$ an algebra. Consider the integral transforms $K_i: U_i (X_i) \longrightarrow V(Y)$ with $i=1,\dots,n+1.$ Then polyconvolution (generalized convolution) of functions $f_1 \in U_1 (X_1), f_2 \in U_2 (X_2),\dots, f_n \in U_n (X_n)$ with weighted $\ga$ for the integral transforms $K_i$ is a multi-linear operators defined by $\textasteriskcentered: \prod\limits_{i=1}^n U_i (X_i)\longrightarrow V(Y),$ denoted by mapping $\underset{\ga}{*} (f_1 , f_2 , \dots, f_n)$ such that the following factorization property is valid
	\begin{equation}\label{1.1}
	K_{n+1} \big[ \underset{\ga}{*} (f_1 , f_2 , \dots, f_n)\big](y)= \ga(y)\prod\limits_{i=1}^n (K_i f_i)(y).
	\end{equation}
	Notice that multiplication on the right-hand side of \eqref{1.1} is understood as multiplication in the algebra $V(Y)$. From \eqref{1.1}, we deduce that
	\begin{equation}\label{1.2}
	\underset{\ga}{*} (f_1 , f_2 , \dots, f_n)(x)=K^{-1}_{n+1} \bigg[\ga(y)\prod\limits_{i=1}^n (K_i f_i)(y) \bigg](x),
	\end{equation}
	where $K^{-1}_{n+1}$ is the inverse transform of $K_{n+1}$.
\end{definition}

Using the Definition \ref{defbyKakichev}, the polyconvolutions generated by various linear operators can be constructed. In particular, this method can be used to construct the polyconvolution involving the Hankel integral transform \cite{Britvina2004,britvina2007} and Fourier-Kontorovich-Lebedev integral transforms \cite{Thao2010Virchenko}, which leads to the study of the operational properties of polyconvolution and applied to solving some classes of integral equations, and systems of integral equations. In this paper, extending the notions in \cite{Kakichev1997}, we introduce a new structure of polyconvolution operator involving the Fourier cosine$(F_c)$-Laplace$(\L)$ integral transforms and apply it to study the solvability in closed-form of classes for integral equations.

The organization of this article is as follows. Besides the introduction, the article has four sections. Section \ref{section2} is devoted to the presentation of the definition of $(F_c, \L)-$polyconvolution. We proved this operator actually satisfies the concept of polyconvolution following Definition \ref{defbyKakichev}. The structurally important properties of this polyconvolution associated with defined spaces are also clearly established. In section \ref{3watson}, we formulated the one-dimensional Watson-type theorem, establish necessary and sufficient conditions for this operator as an isometric isomorphism (unitary) on $L_2 (\R_+)$ space, and get its inverse represented in the conjugate symmetric form. Section 4 consists of two subsections. In Subsection \ref{subsection41}, we prove Young's theorem of the polyconvolution \eqref{eq3.1} and corollaries for estimation. We also prove norm-boundedness in the weighted $L_p (\R_+, \rho_j)$ spaces, which are the contents mentioned in Subsection \ref{subsection42}. Section \ref{apply} deals with applications of the constructed polyconvolution. We study the solvability and unique explicit solution of some classes for integral equations of Toeplitz plus Hankel type and integro-differential equations of Barbashin type by applying the obtained results. Some computational examples can be found in the paper to illustrate the results obtained to ensure their validity and applicability.

\noindent \underline{Background:} We briefly recall some notions used in this article. The Fourier and Fourier cosine transforms of function $f$, denoted by  $(F)$ and  $(F_c)$ respectively, are defined by the integral formulas (see \cite{Sogge1993fourier,bateman1954}) as follows: $(Ff)(y):=\frac{1}{\sqrt{2\pi}}\int\limits_{-\infty}^{+\infty}e^{-ixy}f(x)dx, \ y\in \R$, and $(F_c f)(y):=\sqrt{\frac{2}{\pi}}\int\limits_{0}^{\infty}\cos (xy)f(x)dx,$ $y>0$. The Fourier cosine transform agrees with the Fourier transform if $f(x)$ is an even function. In general, the even part of the Fourier transform of $f(x)$ equals the even part of the Fourier cosine transform of  $f(x)$ in the specified region. The Laplace transform \cite{Sneddon1972} is denoted by $(\L)$ and defined by the integral formula  $(\L f)(y):=\int\limits_{0}^{\infty} e^{-xy} f(x)dx,\ \Re y>0.$
It is well-known that, for all the functions $f \in L_2 (\R)$ then $(Ff)$ is unitary transform on the $L_2 (\R)$ and $\lv Ff\lv_{L_2 (\R)} =\lv f\lv_{L_2 (\R)}$ (refer \cite{Sogge1993fourier,WRudin1987}). The above statement is still true for the $(F_c)$ transform, but the opposite is true for the Laplace to transform $\lv \L f \lv_{L_2 (\R_+)} \leq \sqrt{\pi}\lv f\lv_{L_2 (\R_+)}$ (see \cite{VKTtuan2012}).
In 1941, R.V. Churchill \cite{Churchill1941} considered the convolution of two functions $f$ and $g$ for the Fourier cosine transform defined by 
\begin{equation}\label{tichchapFc}
(f \underset{F_c}{*}g)(x):= \frac{1}{\sqrt{2\pi}} \int\limits_0^\infty f(y)\big[g(x+y)+g(|x-y|) \big]dy,\ x>0.
\end{equation}
This operator is bilinear and commutative (see \cite{MusallamVKt2000,Sneddon1972}). 
If $f,g$  are functions belonging to $L_p (\R_+)$, then  $(f \underset{F_c}{*}g) \in L_p (\R_+)$ \cite{Titchmarsh1986}, with $p=1,2$ and the factorization equality holds
\begin{equation}\label{eq2.2}
F_c (f \underset{F_c}{*}g)(y)=(F_c f)(y)(F_c g)(y), \ y>0.
\end{equation}
Throughout the article, we shall make frequent use of the weighted Lebesgue spaces $L_p (\R_+ , w(x)dx)$, $1\leq p < \infty$ with respect to a positive measure $w(x)dx$ equipped with the norm for which $$\lv f \lv_{L_p (\R_+ , w)}=\bigg( \int\limits_0^{\infty} |f(x)|^p w(x)dx \bigg)^{1/p}<+ \infty.$$ If the weighted function $w=1$ then $L_p (\R_+ ,w)\equiv L_p (\R_+)$. In case $p=\infty$ then the norm of functions defined by $\lv f \lv_{L_{\infty}(\R_+)}= \sup\limits_{x\in \R_+} |f(x)|< +\infty$.

\section{Composition structure of Fouriercosine--Laplace polyconvolution}\label{section2}
In this section, we establish the polyconvolution $\underset{F_c, \L}{*}(f,g,h)(x)$ of three functions $f, g,$ and $h$ constructed for the Fouriercosine $(F_c)$; Laplace $(\L)$ integral transforms. Consider $U_i (X_i)_{i=\overline{1 ,4}}$ on the same field of real numbers and $V(Y)$  as algebraic functions measured in the Lebesgue sense over $\R$. With weight $\ga$ coinciding with constant function 1, transformations $(K_1)=(K_4)=(F_c)$ and $(K_2)=(K_3)=(\L)$, we have definitions.
\begin{definition}\label{dachapFcLaplace}
	The Fourier cosine--Laplace polyconvolution of three functions $f, g, h$ is denoted by $\underset{F_c, \L}{*}(f,g,h)$ and defined by 
	\begin{equation}\label{eq3.1}
	\underset{F_c, \L}{*}(f,g,h)(x):=\frac{1}{\pi}\int\limits_{\R^3_+}\Phi(x,u,v,w) f(u)g(v)h(w)dudvdw,\ x>0.
	\end{equation}
	with kernel function \begin{equation}\label{eq3.2}
	\Phi(x,u,v,w) =\frac{v+w}{(v+w)^2 + (x+u)^2} +\frac{v+w}{(v+w)^2 +(x-u)^2}.
	\end{equation}\end{definition}
\noindent Firstly, we will show that this Definition actually conforms to Kakichev's classical concept of polyconvolution, i.e. operator \eqref{eq3.1} must satisfy the factorization equality in Definition \ref{defbyKakichev} with appropriate integral transformations.
\begin{theorem}\label{theorem31}
	Suppose that $f, g$ and $h$ are arbitrary functions in $L_1 (\R_+)$, then the polyconvolution \eqref{eq3.1}
	for the Fourier cosine and Laplace transforms is well-defined for all $y>0$ as continuous functions and belong to $L_1 (\R_+)$ while keeping the following factorization property holds 
	\begin{equation}\label{eq3.4}
	F_c \big[ \underset{F_c, \L}{*}(f,g,h)\big](y)=(F_c f)(y)(\L g)(y)(\L h)(y),\ \forall y>0.
	\end{equation} Moreover $\underset{F_c, \L}{*}(f,g,h) \in C_0 (\R_+)$, 
	where $C_0(\mathbb{R}_+)$ is the space of bounded continuous functions vanishing at infinity. Besides, the $L_1$-norm estimation of Fourier cosine-Laplace polyconvolution is as follows 
	\begin{equation}\label{eq3.3}
	\lv\underset{F_c, \L}{*}(f,g,h)\lv_{L_1 (\R_+)} \leq \lv f\lv_{L_1 (\R_+)} \lv g\lv_{L_1 (\R_+)} \lv h\lv_{L_1 (\R_+)}. 
	\end{equation}
\end{theorem}
\begin{proof}
	We first prove that $\underset{F_c, \L}{*}(f,g,h)$ belongs to space $L_1 (\R_+)$. From \eqref{eq3.2} we deduce
	\begin{equation}\label{eq3.5}
	\begin{aligned}	
	&\int\limits_{0}^{\infty}|\Phi (x,u,v,w)|dx = \int\limits_{0}^{\infty} \bigg[\frac{v+w}{(v+w)^2 +(x+u)^2}+\frac{v+w}{(v+w)^2+(x-u)^2} \bigg]dx,\ u,v,w>0.\\
	&= \int\limits_{u}^{\infty} \frac{v+w}{(v+w)^2 +t^2}dt+\int\limits_{-u}^{\infty}\frac{u+w}{(v+w)^2 +t^2}dt
	=\int\limits_{-\infty}^{+\infty} \frac{v+w}{(v+w)^2 +t^2}dt = Arc \tan \frac{t}{v+w}\bigg|_{-\infty}^{+\infty} = \pi.
	\end{aligned}
	\end{equation}
	It is easy to see that $\int\limits_{0}^{\infty} |\underset{F_c, \L}{*}(f,g,h)(x)|dx$ is finite. Indeed, from \eqref{eq3.1}, \eqref{eq3.5} and the assumption of the theorem, we obtain
	$$\begin{aligned}
	&\int\limits_{0}^{\infty} |\underset{F_c, \L}{*}(f,g,h)(x) |dx \leq \frac{1}{\pi} \int\limits_{\R^4_{+}} |\Phi (x,u,v,w)|\ |f(u)|\  |g(v)| \ |h(w)|\ dudvdwdx \\
	&=  \bigg(\int\limits_{0}^{\infty}|f(u)|du \bigg)\bigg(\int\limits_{0}^{\infty}|g(v)|dv \bigg)\bigg(\int\limits_{0}^{\infty}|h(w)|dw \bigg)
	= \lv f\lv_{L_1 (\R_+)} \lv g\lv_{L_1 (\R_+)} \lv h\lv_{L_1 (\R_+)} < \infty.
	\end{aligned}$$
	This means that  $\underset{F_c, \L}{*}(f,g,h) \in L_1 (\R_+)$ and estimation \eqref{eq3.3} holds. Next, we prove that the polyconvolution \eqref{eq3.1} satisfying factorization equality \eqref{eq3.4} is valid. 
	Following  formula 2.13.5, page 91 in  \cite{Debnath2006Bhatta}, we obtain $\int\limits_{0}^{\infty}e^{-\a x} \cos(xy)dx=\frac{\a}{\a^2 +y^2},\ \forall \a,y>0.$ 
	Based on this, we can rewrite the equality as follows 
	\begin{equation}\label{eq3.6}
	\begin{aligned}
	\Phi (x,u,v,w)&=\int\limits_{0}^{\infty} e^{-(v+w)y}\{\cos(x+u)y+\cos(x-u)y \}dy
	=2\int\limits_0^{\infty}e^{-(v+w)y} \cos(xy)\ \cos(uy)dy.
	\end{aligned}
	\end{equation}
	Coupling \eqref{eq3.1},\eqref{eq3.6} then
	$\underset{F_c, \L}{*}(f,g,h)(x)= \frac{2}{\pi}\int\limits_{\R^4_+} e^{-(v+w)y} \cos(xy)\ \cos(uy)\ f(u)\ g(v)\ h(w)\ dudvdwdy$. This is equivalent 
	\begin{equation}\label{eq3.7}
	\begin{aligned}
	\underset{F_c, \L}{*}(f,g,h)(x)
	&=	\frac{2}{\pi} \int\limits_{0}^{\infty} \biggl\{\bigg(\int\limits_{0}^{\infty} f(u)\cos(yu)du\bigg) \bigg(\int\limits_{0}^{\infty} g(v)e^{-vy}dv\bigg)\bigg(\int\limits_{0}^{\infty} h(w)e^{-wy}dw\bigg)\biggr\} \cos(xy)dy,\\
	&=\sqrt{\frac{2}{\pi}} \int\limits_{0}^{\infty} (F_c f)(y)(\L g)(y)(\L h)(y)\ \cos(xy)dy,\ \forall x>0.
	\end{aligned}
	\end{equation}	
	By applying the Fourier cosine transformation to both sides of the equality \eqref{eq3.7} and $\underset{F_c, \L}{*}(f,g,h) \in L_1 (\R_+)$, we obtain the factorization property \eqref{eq3.4}. The Riemann-Lebesgue theorem in \cite{Sogge1993fourier}, states that "If $f\in L_1(\R^n)$, then $(Ff)(y) \longrightarrow 0$ as  $|y|\longrightarrow \infty$, and, hence $(Fy )\in C_0 (\R^n)$". This is still true for the  $(F_c)$ transform \cite{Titchmarsh1986}. Therefore, $\underset{F_c, \L}{*}(f,g,h)$ belonging to $C_0 (\R_+)$ is obvious by virtue of Riemann-Lebesgue theorem.
\end{proof}
The following theorems show the  polyconvolution operator \eqref{eq3.1} is always well-defined in different functional spaces.
In the structure of $(F_c , \L)$-polyconvolution, if we fix the functional space of any two functions, namely, let $g,h \in L_1 (\R_+)$ and substitute the assumption $f\in L_1 (\R_+) $ in the theorem \ref{theorem31} by $f \in L_2 (\R_+)$ then we have the following conclusion.
\begin{theorem}\label{theorem32}
	Let $f \in L_2 (\R_+)$ and $g,h$ be arbitrary functions belonging to $L_1 (\R_+)$. Then $\underset{F_c, \L}{*}(f,g,h) \in L_2 (\R_+)$ and for any $x>0$, the following Parseval's  identity holds
	\begin{equation}\label{eq3.8}
	\underset{F_c, \L}{*}(f,g,h)(x)=F_c \big[(F_c f)(y)(\L g)(y)(\L h)y \big](x). 
	\end{equation}
	Besides, the factorization property \eqref{eq3.4} is still valid. 
\end{theorem}
\begin{proof}
	Directly inferred from the formulas \eqref{eq3.5}, \eqref{eq3.6}, for any positive $v,w$,  then $\int\limits_{0}^{\infty} |\Phi(x,u,v,w) |du=\pi$ is finite, and  $|\Phi(x,u,v,w) |\leq 2\int\limits_{0}^{\infty}e^{-(v+w)y}dy<M$, where $M$ is a positive constant. Therefore, by using  H\"older inequality we have the following evaluation
	$$\begin{aligned}
	&|\underset{F_c, \L}{*}(f,g,h)(x) |\\ &\leq \frac{1}{\pi}\biggl\{ \int\limits_{R^3_+} |\Phi(x,u,v,w) |\ |f(u) |^2 |g(v) |\ |h(w) |\ dudvdw \biggr\}^{\frac{1}{2}}\times \biggl\{ \int\limits_{R^3_+} |\Phi(x,u,v,w) |\  |g(v) |\ |h(w) |\ dudvdw \biggr\}^{\frac{1}{2}}\\
	&\leq \sqrt{\frac{M}{\pi}}\biggl\{ \int\limits_{R^3_+}  |f(u) |^2 |g(v) |\ |h(w) |\ dudvdw \biggr\}^{\frac{1}{2}}\times \biggl\{ \int\limits_{R^2_+}   |g(v) |\ |h(w) |\ dvdw \biggr\}^{\frac{1}{2}}\\
	&= \sqrt{\frac{M}{\pi}} \lv f\lv_{L_2 (\R_+)} \lv g\lv_{L_1 (\R_+)} \lv h\lv_{L_1 (\R_+)}<\infty.
	\end{aligned}$$
	\noindent This proves that, the polyconvolution operator \eqref{eq3.1} is existence for any functions $f\in L_2 (\R_+)$ and   $g,h \in L_1 (\R_+)$.  By the same explanations as in \eqref{eq3.7}, we obtain 
	$$\begin{aligned}
	\underset{F_c, \L}{*}(f,g,h)(x)&=\sqrt{\frac{2}{\pi}} \int\limits_{0}^{\infty}(F_c f)(y)(\L g)(y)(\L h)(y)\cos(xy)dy= F_c \big[(F_c f)(y)(\L g)(y)(\L h)(y) \big](x), \ \forall x>0.
	\end{aligned}$$
	Since $f \in L_2 (\R_+)$ then $F_c f$ belongs to the $L_2 (\R_+)$ (refer \cite{Sneddon1972,Titchmarsh1986}). Moreover, $g$ be a function belonging to $L_1 (\R_+)$ then  $|(\L g)(y) |\leq \int\limits_{0}^{\infty} |e^{-vy}|\ | g(v) |dv\leq \int\limits_{0}^{\infty} |g(v)|dv$ is finite $\forall y>0$, this yields $(\L g)(y)$ as a bounded function. Similarly, since $h\in L_1 (\R_+)$ then $ (\L h)(y)$ is bounded. This means that $(F_c f)(y)(\L g)(y)(\L h)(y) \in L_2 (\R_+)$ and $F_c \big[(F_c f)(y)(\L g)(y)(\L h)(y) \big](x)$ belongs to $L_2 (\R_+)$.  Therefore, for all functions $f \in L_2 (\R_+)$ and $g,h \in L_1 (\R_+)$ then $(F_c , \L)$-polyconvolution operator is well-defined in the $L_2 (\R_+)$, and the Parseval's identity \eqref{eq3.8} holds. Furthermore, we know that for the Fourier cosine transform $F_c: L_2 (\R_+)\longleftrightarrow L_2 (\R_+)$ is an isometric isomorphism mapping \cite{Sogge1993fourier}, hence we deduce the factorization equality \eqref{eq3.4} is valid. 
\end{proof}
Before presenting further results, we give a new definition for the functional space class containing $L_1 (\R_+)$   as follows.
\begin{definition}\label{dnAR+}
	Notation $\A(\R_+)$ is the set of all functions defined on $\R_+ $, such that their images via the Laplace transform belong to $ L_{\infty} (\R_ +)$ .
	$$\A(\R_+):=\big\{k(x), x\in \R_+ \ |\ (\L k) \in L_{\infty} (\R_+) \big\}.$$\end{definition}

\noindent It is apparent that with the above definition, for any functions $k$ in $L_1 (\R_+)$, then $k$ also belongs to $\A(\R_+)$.  Opposite, for example consider the function $k(x)= \cos(x)$ or $k(x)=\sin(x)$, then obviously $k(x) \notin L_1 (\R_+)$. However, applying Laplace transform to function $k$, following the formulas 4.7.1, page 150 and 4.7.43, page 154 in \cite{bateman1954} then $(\L \cos(x))(y)=\int\limits_{0}^{\infty} e^{-xy} \cos(x) dx=\frac{y}{1+y^2}$ and $(\L \sin(x))(y)=\int\limits_{0}^{\infty} e^{-xy} \sin(x) dx=\frac{1}{1+y^2}$ are bounded. This implies that $(\L \cos(x)), (\L\sin(x)) \in L_{\infty} (\R_+)$ then $k\in \A(\R_+)$, and therefore we conclude  $L_1 (\R_+) \subseteq \A(\R_+)$.

\begin{theorem}
	If $f$ is a function belonging to $L_2 (\R_+)$ and $g,h \in \A (\R_+)$. Then  $\underset{F_c, \L}{*}(f,g,h)(x)\in L_2 (\R_+)$, the Parseval's  identity \eqref{eq3.8} as well as factorization equality \eqref{eq3.4} still holds.
\end{theorem}
\begin{proof}
	By assuming $f\in L_2 (\R_+)$ then $(F_c f) \in L_2 (\R_+)$ \cite{Titchmarsh1986}. On the other hand $g,h \in \A(\R_+)$, which implies that $(\L g)$ and $(\L h)$ are finite, then $(\L g),(\L h) \in L_{\infty} (\R_+)$. This means that $(\L g)(y), (\L h)(y)$ are bounded functions, and $(F_c f)(y)(\L g)(y)(\L h)(y)$ belongs to $L_2 (\R_+)$, which can be deduced that $$F_c \big[(F_c f)(y)(\L g)(y) (\L h)(y) \big](x) \in L_2 (\R_+).$$ By the same argument as in the Theorem \ref{theorem31} and \ref{theorem32}, we obtain Parseval's equality $\underset{F_c, \L}{*}(f,g,h)(x) = F_c \big[(F_c f)(y)(\L g)(y) (\L h)(y) \big](x) \in L_2 (\R_+),$ thereby using the unitary property of Fouriercosine transform over $ L_2 (\R_+)$, we also get $\eqref{eq3.4}$.
\end{proof}
\section{One-dimensional Watson transform for $(F_c, \L)$ Polyconvolution}\label{3watson}
Following \cite{WatsonGN1933}, G. N. Watson  showed that the Mellin convolution type transforms $g(x)=[K f](x)=\int\limits_{0}^{\infty} k(x y) f(y) d y,$
such that their inverses have the similar form
$f(x)=[\hat{K} g](x)=\int\limits_{0}^{\infty} \hat{k}(x y) g(y) d y,$ where the kernel $ \hat{k}(x)$ is called the conjugate kernel.  Let us consider the above transform  from the point of view of general Fourier transforms (refer \cite{Titchmarsh1986}), it is well-known that the integral transform  $g(x)=\frac{\partial}{\partial x}\int\limits_{0}^{\infty}k_1 (xy)f(y)\frac{dy}{y}$ is an automorphism of the space $L_2 (\R_+)\equiv L_2 (0, \infty)$ and has the symmetric inversion formula if and only if the function $k_1 (x)$ satisfies the conditions

\begin{equation}\label{dkWatson}
\int\limits_{0}^{\infty} k_1 (ax) k_1 (bx) x^{-2}dx = \min (a,b), \hspace{0.3cm} \textup{and} \hspace{0.3cm}
\int\limits_{0}^{\infty} k_1 (ax) \overline{k_1 (bx)}x^{-2}dx =\min (a,b),
\end{equation}
with all positive $a$ and $b$. The function $k_1 (x)$ satisfying the conditions in \eqref{dkWatson} is called one-dimensional Watson kernel \cite{Titchmarsh1986}. It is shown that the conditions \eqref{dkWatson} may be written in
the following equivalent form $k_1 (x)=\frac{x}{2\pi i}\lim\limits_{N\longrightarrow \infty}\int\limits_{\frac{1}{2}-iN}^{\frac{1}{2}+iN}\frac{\Omega(t)}{1-t}x^{-t}dt,$ where the function $\Omega(t)$, defined on the line $\sigma = \{t, \Re t = \frac{1}{2}\}$, fulfills the
conditions $\Omega(t)\Omega(1-t)$ and $|\Omega(t)|=1$. Here \textquotedblleft lim\textquotedblright means the limit with respect to the norm in the space $L_2$. The first equality in condition \eqref{dkWatson} is the classical definition of the Watson kernel \cite{WatsonGN1933}, which leads to a unitary Watson transform when $k_1 (x)$ is a real function. According to this view, we can study the Watson-type integral transform for other multiple convolutions as follows $f\longmapsto g = D(f*k)$, where $D$ is an arbitrary differential operator and $k$ is the known kernel. An example is the previously published works of the Watson-type theorem for convolutions involving the Fourier-cosine \cite{VKT1999JMAA}, Hartley \cite{tuan2022Mediterranean}, Hartley-Fourier \cite{Tuan2022MMA}, and Kontorovich-Lebedev-Fouriersine \cite{tuan2020Ukrainian} transforms. In this section, we study the Watson-type transform for the $(F_c,\L)$ polyconvolution operator by fixing a function, and letting the remaining functions vary in certain functional spaces. We establish  necessary and sufficient conditions for this operator to be unitary on $L_2 (\R_+)$ and get its inverse represented in the conjugate symmetric form. Consider the operator

$$\begin{aligned}
\T_{\eta, \xi} : L_2 (\R_+)&\longrightarrow L_2 (\R_+)\\
f&\longmapsto \vp= (\T_{\eta, \xi} f) =D( \underset{F_c, \L}{*}(f, \eta,\xi)),
\end{aligned}$$
\noindent where $D$ is the second-order differential operator defined by $D:=\left(I-\frac{d^2}{dx^2}\right)$ and image
\begin{equation}\label{eq4.1}
\begin{aligned}
\vp(x)=(\T_{\eta ,\xi}f)(x)= \left(1-\frac{d^2}{dx^2}\right)\frac{1}{\pi} \int\limits_{\R^3_+} \Phi(x,u,v,w)\ f(u)\ \eta(v)\xi(w)\ dudvdw,\ \forall x>0.
\end{aligned}
\end{equation}
Here $\Phi(x,u,v,w)$ is defined by \eqref{eq3.2} and $\eta, \xi$ are given functions.
\begin{theorem}\label{theorem4.1}
	Suppose that $\eta , \xi \in \A(\R_+)$ are given functions and satisfy 
	\begin{equation}\label{dieukien4.2}
	|(\L \eta)(y) (\L \xi)(y) |=\frac{1}{1+y^2},\ \forall y>0.
	\end{equation}
	Then the condition \eqref{dieukien4.2} is the necessary and sufficient condition for operator $(\T_{\eta, \xi})$ to become an isometric isomorphism mapping on $L_2 (R_+)$. Moreover, the inverse operator of $(\T_{\eta, \xi})$ has a symmetric form and is represented by 
	$$\begin{aligned}
	f(x)=(\T_{\bar\eta, \bar\xi} \vp)(x)&=D(\underset{F_c, \L}{*}(f, \bar\eta,\bar\xi))(x)\\
	&= \left(1-\frac{d^2}{dx^2}\right)\frac{1}{\pi} \int\limits_{\R^3_+} \Phi(x,u,v,w)\ \vp(u)\ \bar\eta(v)\bar\xi(w)\ dudvdw,\ \forall x>0,
	\end{aligned}$$
	where $\bar\eta,\bar\xi$ are complex conjugate functions of $\eta, \xi$ respectively.
\end{theorem}
In order to prove Theorem \ref{theorem4.1}, we need the following auxiliary lemma.

\begin{lemma}\label{bode4.1}
	If $\eta, \xi$ are functions belonging to  $\A (\R_+)$ and satisfy 
	\begin{equation}\label{dieukien4.3}
	(1+y^2) |(\L \eta)(y)(\L\xi)(y) | \ \textup{  bounded for any}\ y>0.
	\end{equation}
	Then, for any function $f$ belonging to the $L_2 (\R_+)$, we obtain
	\begin{equation}\label{eq4.4}
	D(\underset{F_c, \L}{*}(f,g,h))(x)= F_c \big[(1+y^2)(F_c f)(y)(\L \eta)(y) (\L \xi)(y) \big](x) \in L_2 (\R_+),	
	\end{equation} where $D$ is the second order differential operator having the form $D=\left(I-\frac{d^2}{dx^2}\right)$.
\end{lemma}
\begin{proof}
	One obvious thing is that, the functions $f(y), yf(y), y^2 f(y), \dots, y^n f(y)$ belong to
	the space $L_2 (\R)$ if and only if $(Ff)(x),\frac{d}{dx} (Ff)(x), \frac{d^2}{dx^2} (Ff)(x), \dots, \frac{d^n}{dx^n} (Ff)(x)\in$  $L_2 (\R)$ (refer Theorem 68, page92, in \cite{Titchmarsh1986}). This confirmation is still true for the Fourier cosine $(F_c)$ transform on $L_2 (\R_+)$. Furthermore, we have	
	$\frac{d}{dx} (F_c f)(x)=\sqrt{\frac{2}{\pi}} \int\limits_{0}^{\infty} f(y)\left(\frac{d}{dx}\cos(xy)\right)dy=\sqrt{\frac{\pi}{2}}\int\limits_{0}^{\infty} (-y)f(y)\sin(xy)dy,$ and
	$$\frac{d^2}{dx^2}(F_c y)(x)=\sqrt{\frac{\pi}{2}}\int\limits_{0}^{\infty}(-y)f(y)\left(\frac{d}{dx}\sin(xy)\right)dy=F_c (-y^2f(y))(x) \in L_2 (\R_+).$$
	That means, if $f (y), y^2 f (y)$ belong to $L_2 (R_+)$, then we obtain \begin{equation}\label{eq4.5}
	\left(1-\frac{d^2}{dx^2}\right)\big(F_c f(y)\big)(x)=F_c \big[(1+y^2)f(y)\big](x).\end{equation}
	By condition \eqref{dieukien4.3}, we deduce that $(1+y^2)(\L \eta)(y)(\L \xi)(y)$  is a bounded function $\forall y>0$, and $f\in L_2 (\R_+)$ implies that $F_c f\in L_2 (\R_+)$. Thus, $(1+y^2)(F_c f)(y)(\L \eta)(y)(\L \xi)(y)$ belongs to the $L_2(R_+)$. Therefore, deducing directly from Parseval's identity \eqref{eq3.8} and \eqref{eq4.5}, we obtain
	
	$$\begin{aligned}
	D( \underset{F_c, \L}{*}(f, g,h))(x)&= \left(1-\frac{d^2}{dx^2}\right) \big\{F_c\big[(F_c f)(y)(\L \eta)(y)(\L \xi)(y) \big](x) \big\}\\
	&= F_c \big[(1+y^2)(F_c f)(y) (\L \eta)(y)(\L \xi)(y)\big](x) \in L_2 (\R_+).
	\end{aligned}$$
\end{proof}

\begin{proof}\textit{of Theorem \ref{theorem4.1}}\\
	\textit{Necessary condition:} From \eqref{dieukien4.2}, it can be seen that $(1+y^2)(\L \eta )(y)(\L \xi)(y)$ is bounded, together with $F_c f$ belonging to $L_2 (\R_+)$, we deduce that $(1+y^2)(F_c f)(y)(\L \eta)(y)(\L \xi)(y) \in L_2 (\R_+)$. Applying the equality \eqref{eq4.4} to \eqref{eq4.1}, we have
	\begin{equation}\label{eq4.6}
	\begin{aligned}
	\vp(x)=(\T_{\eta,\xi}f)(x)
	&=D (\underset{F_c, \L}{*}(f,\eta,\xi))(x)\\
	&= F_c \big[(1+y^2)(F_c f)(y)(\L \eta)(y)(\L \xi)(y) \big](x) \in L_2 (\R_+).
	\end{aligned}
	\end{equation}
	According to Theorem 9.13 in \cite{WRudin1987}, for any $f$ be a function belonging to the $L_2 (\R)$, then $\lv Ff\lv_{L_2 (\R)}=\lv f\lv_{L_2 (\R)}$.  This is still true for the $(F_c)$ transform which implies that $\lv F_c f\lv_{L_2 (\R_+)}=\lv f\lv_{L_2 (\R_+)}, \forall f\in L_2 (\R_+)$.
	Applying this to equality \eqref{eq4.6} and under the condition \eqref{dieukien4.2}, we obtain
	$$\begin{aligned}
	\lv\vp \lv_{L_2 (\R_+)}= \lv\T_{\eta ,\xi}f \lv_{L_2 (\R_+)}&= \lv D (\underset{F_c, \L}{*}(f,\eta,\xi)) \lv_{L_2 (\R_+)}\\
	&=(1+y^2)|(\L \eta)(y) (\L \xi)(y) |\ \lv F_c f\lv_{L_2 (\R_+)} =\lv F_c f \lv_{L_2 (\R_+)} =\lv f\lv_{L_2 (\R_+)}.
	\end{aligned}$$
	This means $\T_{\eta,\xi}$ is an isometric isomorphism (unitary) transformation  in $L_2(R_+)$. Now, we need to show the inverse operator of $(\T_{\eta, \xi})$ has a symmetric form. Indeed, in $L_2$, by apply the unitary property of mapping $F_c: L_2 (\R_+)\longleftrightarrow L_2 (\R_+)$ (see \cite{Titchmarsh1986}) for equality \eqref{eq4.6}, we get
	\begin{equation}\label{eq4.7}
	(F_c \vp)(y)=F_c (\T_{\eta ,\xi}f)(y)=(1+y^2)(F_c f)(y)(\L \eta)(y)(\L\xi)(y) \in L_2 (\R_+).
	\end{equation}
	On the other hand, since $(\overline{\L \eta})(y)=(\L \bar\eta)(y)$ and $(\overline{\L \xi})(y)=(\L \bar\xi)(y)$, we have  $$(1+y^2)| (\overline{\L \eta})(y)(\overline{\L \xi})(y)|=(1+y^2)(\L \bar\eta)(y)(\L \bar\xi)(y).$$ Due to \eqref{dieukien4.2}, it follows that $(1+y^2)| (\overline{\L \eta})(y)(\overline{\L \xi})(y)|=1$. This leads to $(\overline{\L \eta})(y)(\overline{\L \xi})(y)(F_c \vp)(y)=\frac{1}{1+y^2}(F_c f)(y)$ based on \eqref{eq4.7}, equivalent to $(\L \bar\eta)(y)(\L \bar\xi)(y)(F_c \vp)(y)=\frac{1}{1+y^2} (F_c f)(y)$. We obtain
	\begin{equation}\label{eq4.8}
	(1+y^2)(\L \bar \eta)(y)(\L \bar\xi)(y)(F_c \vp)(y)=(F_c f)(y).
	\end{equation}
	\noindent By the same reasoning as above, we get $F_c \vp \in L_2 (\R_+)$ and $(1+y^2)(\L \bar \eta)(y)(\ L \bar \xi)(y)$ is a bounded function, implying that $(1+y^2)(\L \bar \eta)(y)(\L \bar \xi)(y)(F_c \vp) (y) \in L_2 (\R_+)$. Combining \eqref{eq4.8} and \eqref{eq4.4}, we obtain
	$$\begin{aligned}
	f(x)&=F_c \big[(1+y^2)(\L \bar\eta)(y) (\L \bar \xi)(y)(F_c \vp)(y) \big](x)
	=D(\underset{F_c,\L}{*}(\vp, \bar \eta, \bar \xi) )(x)=(\T_{\bar \eta, \bar \xi}\vp)(x).
	\end{aligned}$$
	\noindent \textit{Sufficient condition:} Assume that $\T_{ \eta,  \xi}$ is a unitary operator on $L_2(\R_+)$ and has the inverse operator $\T_{\bar \eta, \bar \xi}$. We need to show that the functions $\eta, \xi$ must satisfy the condition \eqref{dieukien4.2}. Indeed, since $\T_{ \eta,  \xi}$ has the unitary property on $L_2(\R_+)$, then for any functions $\vp$ belonging to $L_2 (\R_+)$, we obtain
	$$\begin{aligned}
	\lv F_c \vp	\lv_{L_2 (\R_+)} =\lv F_c(\T_{ \eta,  \xi} f) \lv_{L_2 (\R_+)}&=	\lv  \vp\lv_{L_2 (\R_+)}
	=\lv F_c\big[(1+y^2)(F_c f)(y)(\L \eta)(y)(\L \xi)(y) \big]	\lv_{L_2 (\R_+)}\\
	&=(1+y^2)|(\L \eta )(y)(\L \xi)(y) |.\ \lv F_c f\lv_{L_2 (\R_+)}\\&=\lv F_c f\lv_{L_2 (\R_+)}=\lv f \lv_{L_2 (\R_+)},\ \forall f \in L_2 (\R_+).
	\end{aligned}$$
	This shows that there exists a multiplication operator of the form $\M_{\Theta}[.]$  defined by  $\M_{\Theta}[f](y):=\Theta(y).f(y)$, where the function $\Theta(y)=(1+y^2).|(\L \eta )(y)(\L \xi)(y)|,\ \forall y>0$. The above expression can be rewritten as $\lv F_c f\lv_{L_2 (\R_+)}=\lv\M_{\Theta} (F_c f) \lv_{L_2 (\R_+)}$ for any $f \in L_2 (\R_+)$. This  means that $\M_{\Theta}[.]$ is an isometric isomorphism (unitary) on the $L_2 (\R_+)$, and this happens if and only if $|(\L \eta )(y)(\L \xi)(y)| =\frac{1}{1+y^2}.$ In conclusion, both $\eta$ and $\xi$ must satisfy the condition \eqref{dieukien4.2}.\end{proof}

\begin{remark}
	It should be emphasized that condition \eqref{dieukien4.2} is indeed narrower than condition \eqref{dieukien4.3}, and these conditions are well-defined with the given assumptions. To make it clear, we consider an example of the pair of functions  $\eta, \xi$ that satisfies condition \eqref{dieukien4.3} but does not satisfy the condition of \eqref{dieukien4.2}. Let $\eta= i\sin t$ and $\xi =\cos t$. It is easy to check these functions do not belong to $L_1 (\R_+)$, according to the formulas 4.7.1, page 150 and 4.7.43, page 154 in \cite{bateman1954}, we obtain both $(\L i\sin t)(y)=\frac{i}{1+y^2}$ and $(\L \cos t)(y)=\frac{y}{1+y^2}$ respectively finite for any $y>0$. This means that $(\L i\sin t)$ and $(\L \cos t) \in L_{\infty}(\R_+)$, implying $\eta, \xi$ belong to the $\A (\R_+)$. Therefore, $(1+y^2)|(\L i \sin t)(y)(\L \cos t)(y)|=\frac{y}{(y^2+1)^2}$ is bounded. In this case, we conclude that the pair $ \eta, \xi $ satisfies \eqref{dieukien4.3} but  the condition \eqref {dieukien4.2} is invalid.
	
	A straightforward instance shows that the pair $\eta, \xi$  actually satisfies the condition \eqref{dieukien4.2}. Let $\eta =e^{it}$ and $\xi =e^{-it}$.  According to the formula 4.5.1, page 143 in \cite{bateman1954}, then $(\L e^{it})(y)=\frac{1}{y-i}$ is finite and $ (\L e^{-it})(y)=\frac{1}{y+i}$ is finite, implying that $(\L e^{\pm it}) \in L_{\infty} (\R_+)$. Therefore $\eta, \xi \in \A (\R_+)$ and $(\L e^{it})(y)(\L e^{-it})(y) =\frac{1}{1+y^2}, \forall y>0$, deduce these functions indeed satisfy the condition \eqref{dieukien4.2}.
	The condition \ref{dieukien4.2} would not exist if we substituted the assumptions of Theorem \ref{theorem4.1} by $\eta,\xi \in L_1 (\R_+)$.
\end{remark}
\begin{remark}
	In the general case, Theorem \ref{theorem4.1} can be studied by replacing the second-order differential operator by $D$ with an arbitrary differential operator of order $2n$ with coefficients in the following form  
	$D^{2n}:=\sum\limits_{k=0}^{n}(-1)^k a_k \frac{d^{2k}}{dx^{2k}},\ n \in \N.$
	Then, the condition \eqref{dieukien4.2} becomes $| (\L \eta)(y)(\L \xi)(y)|=\frac{1}{\P_{2n}(y)},$ which is the necessary and sufficient condition for operator $(\T_{\eta, \xi})$ to be unitary on $L_2 (R_+)$. Here,  $\P_{2n}(y)=\sum\limits_{k=0}^{n} a_k y^{2k}$ is a polynomial with real coefficients without real zero-points. It can be seen that, the second-order differential operator $D$ becomes a special case with $n=1, a_0 =a_1 =1.$
\end{remark}
\section{Young-type theorem and estimation in  $L_p (\R_+ , \rho_j)$ weighted  spaces}
For the Fourier convolution 
$(f\underset{F}{*} g)(x):=\frac{1}{\sqrt{2 \pi}}\int\limits_{-\infty}^{\infty} f(x-y)\ g(y)dy,\ x\in \R$,
W.H. Young \cite{YoungWH1912} obtained the following inequality 
$\lv f \underset{F}{*} g  \lv_{L_r (\R)} \leq \lv f\lv_{L_p (\R)} \lv g\lv_{L_q (\R)},$
where $p, q, r > 1$  such that $\frac{1}{p} +\frac{1}{q}=1+\frac{1}{r}$ and $f \in L_p (\R),\ g\in L_q(\R)$. Afterwards, R.A. Adams  and J.F. Fournier generalized Young's inequality for the Fourier convolution (\cite{AdamsFournier2003sobolev}, Theorem 2.24, page 33) to include a weight
$\big|\int\limits_{\R^n} (f\underset{F}{*} g)(x). w(x)dx \big| \leq \lv f\lv_{L_p (\R^n)} \lv g\lv_{L_q (\R^n)} \lv w\lv_{L_r (\R^n)}$
where  $p, q, r>1$ such that $\frac{1}{p} +\frac{1}{q}+\frac{1}{r}=2$ with $f \in L_p (\R^n), g\in L_q (\R^n),$ and $w \in L_r (\R^n)$.  Notice that for the important case $f, g \in L_2 (\R)$,  Young's inequality does not hold.  In \cite{Saitoh2000}, Saitoh derived a weighted $L_{p}(\mathbb{R},|\rho_j|)$ norm inequality for the Fourier convolution of the following form:
$$\big\lv\left( (F_1 \rho_1) \underset{F}{*}(F_2 \rho_2)\right) (\rho_1 \underset{F}{*} \rho_2)^{\frac{1}{p}-1}\big\lv_{L_p (\R)} \leq \lv F_1 \lv_{L_p (\R,|\rho_1|)} \lv F_2\lv_{L_p (\R,|\rho_2|)},\ \textup{with}\ p>1,$$ where $\rho_j $ are non-vanishing functions, $F_j \in L_p (\R, |\rho_j|),\, j=1,2$. Here, the norm of $F_j$ in the weighted space $L_{p}(\mathbb{R}, \rho_j)$ is understood as $
\lv F_j  \lv_{L_{p}(\mathbb{R}, \rho_j)}=\big\{\int\limits_{-\infty}^{\infty}|F_j(x)|^{p} \rho_j(x) \mathrm{d} x\big\}^{\frac{1}{p}}.$
The reverse weighted $L_p$-norm inequality for Fourier convolution has also been investigated. Unlike Young's inequality, Saitoh's inequality also holds in case $p=2$, which is the most obvious difference between these two inequalities. Furthermore, in many cases of interest, the convolution is given in the form $
\rho_{2}(x) \equiv 1,$\ $F_{2}(x)=G(x),
$
where $G(x-\xi)$ is some Green's functions. Then the above inequality becomes
$
\|(F \rho) \underset{F}{*} G\|_{L_{p}(\mathbb{R})} \leqslant\|\rho\|_{L_{1}\left(\mathbb{R}_{+}\right)}^{1-1/p}\|G\|_{L_{p}(\mathbb{R})}\|F\|_{L_{p}(\mathbb{R},|\rho|)},
$
where $\rho, F$ and $G$ are such that the right-hand side  is finite.  Saitoh's inequality can be applied to estimating  the solution to a parabolic integro-differential equation \cite{hoangtuan2017thaovkt} modeling a scattered acoustic field. Based on the above aspects, we will obtain certain norm inequalities for polyconvolution \eqref{eq3.1} in a very general framework, and estimation in $L_p$ weighted space.  Some techniques used in this proof of our theorem come from \cite{tuan2020Ukrainian,TuanVKTuan2023ITSF}. We follow closely the strategy of these results, especially how to choose the bounded linear functional in \cite{TuanVKTuan2023ITSF} to apply Riesz's representation theorem skillfully.
\subsection{Young-type theorem for $(F_c, \L)-$polyconvolution operator}\label{subsection41}
\begin{theorem}\label{Young-typeTheorem}
	Let $p,q,r$, and $s$ be real numbers in open interval $(1,\infty)$ such that $1\textfractionsolidus p + 1\textfractionsolidus q + 1\textfractionsolidus r+ 1\textfractionsolidus s=3$.	
	For any functions $f\in L_p (\R), g\in L_q (\R_+, (w+x)^{q-1}), h\in L_r (\R_+, (v+x)^{r-1}),$ and $k\in L_s (\R_+)$, the following inequality holds true for
	\begin{equation}\label{eq5.1}
	\bigg|\int\limits_{0}^{\infty} (\underset{F_c ,\L}{*}(f,g,h))(x). k(x)dx \bigg|\leq w^{(1-q)\textfractionsolidus q} v^{(1-r)\textfractionsolidus r} \lv f\lv_{L_p (\R_+)}\lv g\lv_{L_q (\R_+, (w+x)^{q-1})} \lv h\lv_{L_r (\R_+, (v+x)^{r-1})}\lv k \lv_{L_s (\R_+)},
	\end{equation}	
\end{theorem}
\begin{proof}
	Based on \eqref{eq3.2}, $|\Phi(x,u,v,w)|\leq \frac{2}{w+v}$. We have
	$
	\int\limits_{0}^{\infty} \big| \frac{\Phi(x,u,v,w)}{v+w}\big|dv\leq 2 \int\limits_{0}^{\infty} \frac{1}{(v+w)^2}dv\leq 2\int\limits_{0}^{\infty} \frac{1}{v^2 + w^2}dv =\frac{2}{w} Arc\tan \frac{v}{w}\big|_{0}^{\infty}=\frac{\pi}{w}, \forall w>0.
	$
	On a similar way, we also obtain $\int\limits_{0}^{\infty} \big| \frac{\Phi(x,u,v,w)}{v+w}\big|dw \leq \frac{\pi}{v}$  finite $\forall v>0.$ 
	Let  $p_1,q_1,r_1,s_1$  be the conjugate exponentials of $p,q,r,s$ respectively. This means that  $\frac{1}{p}+\frac{1}{p_1}=\frac{1}{q}+\frac{1}{q_1}=\frac{1}{r}+\frac{1}{r_1}+\frac{1}{s}+\frac{1}{s_1}=1$ and together with the assumption of theorem, we get the correlation between exponential numbers as follows
	\begin{align}\label{eq5.3}
	\begin{cases}
	&1\textfractionsolidus p_1 + 1\textfractionsolidus q_1 + 1\textfractionsolidus r_1 + 1\textfractionsolidus s_1=1 \\
	&p\left(\frac{1}{q_1}+\frac{1}{r_1}+\frac{1}{s_1}\right)=q\left(\frac{1}{p_1}+\frac{1}{r_1}+\frac{1}{s_1}\right)=r\left(\frac{1}{p_1}+\frac{1}{q_1}+\frac{1}{s_1}\right)=s\left(\frac{1}{p_1}+\frac{1}{q_1}+\frac{1}{r_1}\right)=1.\\
	&(q-1)\left(\frac{1}{p_1}+\frac{1}{r_1}+\frac{1}{s_1}\right)-\frac{1}{q_1} =(r-1)\left(\frac{1}{p_1}+\frac{1}{q_1}+\frac{1}{s_1}\right)-\frac{1}{r_1}=0.
	\end{cases}
	\end{align}
	
	\noindent  For simplicity, we denote $\Omega = \R^4_+$ and without loss of generality, for all $(x,u,v,w) \in \Omega$, we put
	$$\begin{aligned}
	T_1 (x,u,v,w) &= |g(v)|^{\frac{q}{p_1}}\ |h(w)|^{\frac{r}{p_1}} |k(x)|^{\frac{s}{p_1}}\ |v+w|^{\frac{(q-1)+(r-1)}{p_1}}\ |\Phi (x,u,v,w)|^{\frac{1}{p_1}},\\
	T_2 (x,u,v,w) &= |f(u)|^{\frac{p}{q_1}}\ |h(w)|^{\frac{r}{q_1}} |k(x)|^{\frac{s}{q_1}}\ |v+w|^{\frac{r-1}{q_1}}\ \big|\frac{\Phi (x,u,v,w)}{v+w}\big|^{\frac{1}{q_1}},\\
	T_3 (x,u,v,w) &= |f(u)|^{\frac{p}{r_1}}\ |g(v)|^{\frac{q}{r_1}} |k(x)|^{\frac{s}{r_1}}\ |v+w|^{\frac{q-1}{r_1}}\ \big|\frac{\Phi (x,u,v,w)}{v+w}\big|^{\frac{1}{r_1}},\\
	T_4 (x,u,v,w) &= |f(u)|^{\frac{p}{s_1}}\ |g(v)|^{\frac{q}{s_1}} |k(x)|^{\frac{r}{s_1}}\ |v+w|^{\frac{(q-1)+(r-1)}{s_1}}\ |\Phi (x,u,v,w)|^{\frac{1}{s_1}}.
	\end{aligned}$$
	From \eqref{eq5.3}, we deduce that
	
	\begin{equation}\label{eq5.4}
	\prod\limits_{i=1}^4 T_i (x,u,v,w) = |f(u)|\ |g(v)|\ |h(w)|\ |k(x)|\ |\Phi (x,u,v,w)|,\ \forall (x,u,v,w) \in \Omega.
	\end{equation}
	Based on the assumption of $g\in L_q (\R_+, (w+x)^{q-1}), h\in L_r (\R_+, (v+x)^{r-1}), k\in L_s (\R_+)$, using Fubini's theorem,  we obtain  $L_{p_1} (\Omega)$-norm estimation for the operator $T_1$ as follows
	\begin{equation}\label{eq5.5}
	\begin{aligned}
	\lv T_1 \lv^{p_1}_{L_{p_1}(\Omega)} &\leq \int\limits_\Omega |g(v)|^q |h(w)|^r |k(x)|^s |v+w|^{(q-1)} |v+w|^{(r-1)} |\Phi (x,u,v,w)|\ dudvdwdx\\
	&\leq \pi \bigg(\int\limits_{0}^{\infty} |g(v)|^q |v+w|^{(q-1)} dv\bigg) \bigg(\int\limits_{0}^{\infty} |h(w)|^r |v+w|^{(r-1)} dw \bigg) \bigg(\int\limits_{0}^{\infty} |k(x)|^s dx \bigg)\\&
	= \pi \lv g \lv^q_{L_q (\R_+,(x+w)^{q-1})} \lv h \lv^r_{L_r (\R_+,(x+v)^{r-1})} \lv k\lv^s_{L_s (\R_+)}.
	\end{aligned}
	\end{equation}
	Similar to what we did with the evaluation   of $T_1$, we also get the norm estimation of $T_4$ on $L_{s_1} (\Omega)$ as follows
	\begin{equation}\label{eq5.6}
	\lv T_4 \lv^{s_1}_{L_{s_1}(\Omega)} \leq \pi \lv f\lv^p_{L_p (\R_+)} \lv g\lv^q_{L_q (\R_+,(x+w)^{q-1})} \lv h\lv^r_{L_r (\R_+,(x+v)^{r-1})}.
	\end{equation}
	For the evaluation of operator $T_2$ on $L_{q_1} (\Omega)$, since $\int\limits_{0}^{\infty} \big| \frac{\Phi(x,u,v,w)}{v+w}\big|dv \leq \frac{\pi}{w}$, using Fubini's theorem  we have 
	\begin{equation}\label{eq5.7}\begin{aligned}
	\lv T_2 \lv^{q_1}_{L_{q_1}(\Omega)} &=\int\limits_{\Omega} |f(u)|^p |h(w)|^r |k(x)|^s |v+w|^{r-1} \bigg|\frac{\Phi(x,u,v,w)}{v+w} \bigg|dudvdwdx\\
	&\leq\frac{\pi}{w}\bigg(\int\limits_{0}^{\infty} |f(u)|^p du \bigg) \bigg(\int\limits_{0}^{\infty} |h(w)|^r |v+w|^{r-1} dw \bigg) \bigg(\int\limits_{0}^{\infty} |k(x)|^s dx \bigg)\\&
	=\frac{\pi}{w}\lv f\lv^p_{L_p (\R_+)} \lv h\lv^r_{L_r (\R_+,(v+x)^{r-1})} \lv k\lv^s_{L_s (\R_+)},\ w>0.
	\end{aligned}
	\end{equation}
	And $L_{r_1} (\Omega)$-norm estimation for the operator $T_3$ has the following form
	\begin{equation}\label{eq5.8}
	\lv T_3 \lv^{r_1}_{L_{r_1}(\Omega)} \leq \frac{\pi }{v}\lv f\lv^p_{L_p (\R_+)} \lv g\lv^q_{L_q (\R_+,(x+w)^{q-1})} \lv k\lv^s_{L_s (\R_+)}.
	\end{equation}
	Combining \eqref{eq5.5}\eqref{eq5.6}\eqref{eq5.7} and  \eqref{eq5.8}, under condition \eqref{eq5.3}, we obtain
	
	\begin{equation}\label{eq5.9}
	\begin{aligned}
	&\lv T_1 \lv_{L_{p_1}(\Omega)}\lv T_2 \lv_{L_{q_1}(\Omega)}\lv T_3 \lv_{L_{r_1}(\Omega)}\lv T_4 \lv_{L_{s_1}(\Omega)}\\ &\leq \pi^{\left(\frac{1}{p_1}+\frac{1}{q_1}+\frac{1}{r_1}+\frac{1}{s_1}\right)} \left(\frac{1}{w}\right)^{\frac{1}{q_1}}\left(\frac{1}{v}\right)^{\frac{1}{r_1}} \lv f\lv^{p\left(\frac{1}{q_1}+\frac{1}{r_1}+\frac{1}{s_1}\right)}_{L_p (\R_+)} 
	\lv g\lv^{q\left(\frac{1}{p_1}+\frac{1}{r_1}+\frac{1}{s_1}\right)}_{L_q (\R_+,(x+w)^{q-1})}
	\lv h\lv^{r\left(\frac{1}{p_1}+\frac{1}{q_1}+\frac{1}{s_1}\right)}_{L_q (\R_+,(x+v)^{r-1})}
	\lv k \lv^{s\left(\frac{1}{p_1}+\frac{1}{q_1}+\frac{1}{r_1}\right)}_{L_s (\R_+)}\\
	&=\pi  \left(\frac{1}{w}\right)^{\frac{1}{q_1}} \left(\frac{1}{v}\right)^{\frac{1}{r_1}} \lv f\lv_{L_q (\R_+)} \lv g\lv_{L_q (\R_+,(x+w)^{q-1})} \lv h\lv_{L_r (\R_+(x+v)^{r-1})} \lv k\lv_{L_s (\R_+)}.
	\end{aligned}
	\end{equation}
	\noindent Moreover, from \eqref{eq3.1} and  \eqref{eq5.4}, we get $$\begin{aligned}
	\big|\int\limits_0^{\infty} (\underset{F_c , \L}{*}(f,g,h))(x).k(x)\big|&\leq \frac{1}{\pi} \int\limits_{\R^4_+} \big|\Phi(x,u,v,w)\big|\ |f(u)|\ |g(v)|\ |h(w)|\ |k(x)|\ dudvdwdx\\
	&=\frac{1}{\pi} \int\limits_{\R^4_+}  \prod\limits_{i=1}^4 T_i (x,u,v,w) \ dudvdwdx.\end{aligned}$$ Since $p_1, q_1, r_1, s_1$ are the conjugate exponentials then $\frac{1}{p_1}+\frac{1}{q_1}+\frac{1}{r_1}+\frac{1}{s_1}=1$, applying Hölder's inequality to four operators and following \eqref{eq5.9}, we obtain
	$$\begin{aligned}
	&\bigg|\int\limits_0^{\infty} (\underset{F_c , \L}{*}(f,g,h))(x).k(x)\bigg|= \frac{1}{\pi}\int\limits_{\R^4_+}  \prod\limits_{i=1}^4 T_i (x,u,v,w) \ dudvdwdx\\
	&\leq \frac{1}{\pi}\biggl\{ |T_1|^{p_1 }\ dudvdwdx\biggr\}^{\frac{1}{p_1}}	\times \biggl\{ |T_2 |^{q_1}\ dudvdwdx\biggr\}^{\frac{1}{q_1}}	
	\times \biggl\{ |T_3|^{r_1 }\ dudvdwdx\biggr\}^{\frac{1}{r_1}}	\times \biggl\{ |T_4|^{s_1 }\ dudvdwdx\biggr\}^{\frac{1}{s_1}}\\
	&=\frac{1}{\pi} \lv	T_1 \lv_{L_{p_1}(\Omega)}\lv	T_2 \lv_{L_{q_1}(\Omega)} \lv	T_3 \lv_{L_{r_1}(\Omega)} \lv	T_4 \lv_{L_{s_1}(\Omega)}\\
	&\leq \left(\frac{1}{w}\right)^{\frac{1}{q_1}} \left(\frac{1}{v}\right)^{\frac{1}{r_1}} \lv f \lv_{L_p (\R_+)} \lv g \lv_{L_q (\R_+,(x+w)^{q-1})} \lv h \lv_{L_r (\R_+,(x+v)^{r-1})} \lv k \lv_{L_s (\R_+)}.
	\end{aligned}$$
	\noindent Finally, since $1\textfractionsolidus q + 1\textfractionsolidus q_1=1$, then  $\left(\frac{1}{w}\right)^{\frac{1}{q_1}}=w^{(1-q)\textfractionsolidus q}$, and  $1\textfractionsolidus r + 1\textfractionsolidus r_1=1$, then $\left(\frac{1}{v}\right)^{\frac{1}{r_1}}=v^{(1-r)\textfractionsolidus r}$,  getting the desired conclusion.
\end{proof}
In case the given function $k(x)$ becomes $(F_c,\L)$-polyconvolution operator \eqref{eq3.1}, then the following Young type inequality is a direct consequence of Theorem \ref{Young-typeTheorem}.
\begin{corollary}\label{hequa5.1}
	Let $p,q,r,s > 1$ be real numbers, satisfy  $1\textfractionsolidus p+1\textfractionsolidus q+1\textfractionsolidus r=1\textfractionsolidus s+2$. If $f \in L_p (\R_+),g\in L_q (\R_+,(w+x)^{q-1})$ and $h\in L_r (\R_+, (v+x)^{q-1})$, then the polyconvolution $\underset{F_c , \L}{*}(f,g,h)$ is well-defined and belongs to $L_s (\R_+)$.  Moreover, the following inequality holds
	\begin{equation}\label{eq5.10}
	\lv \underset{F_c , \L}{*}(f,g,h)\lv_{L_s (\R_+)} \leq w^{(1-q)\textfractionsolidus q} v^{(1-r)\textfractionsolidus r} \lv f\lv_{L_p (\R_+)} \lv g\lv_{L_q (\R_+,(w+x)^{q-1}}\lv h\lv_{L_r (\R_+,(v+x)^{r-1}}.
	\end{equation}
\end{corollary}
\begin{proof}
	Let $s_1$ be the conjugate exponent of $s$, i.e $\frac{1}{s}+\frac{1}{s_1} =1$. From the assumptions of Corollary \ref{hequa5.1}, we have $\frac{1}{p}+\frac{1}{q}+\frac{1}{r}+\frac{1}{s_1}=3$,  which shows the numbers $p,q,r,$ and $s_1$  satisfy the conditions of  Theorem \ref{Young-typeTheorem} (with $s$ being replaced by $s_1$). Therefore, if $f \in L_p (\R_+),g\in L_q (\R_+,(w+x)^{q-1})$ and $h\in L_r (\R_+, (v+x)^{q-1})$, then the linear functional	$$Lk := \int\limits_{0}^{\infty} \big(\underset{F_c , \L}{*}(f,g,h)\big)(x)\ .k(x)dx $$ is bounded in $L_{s_1}(\R_+)$. Consequently, by the Riesz's representation theorem \cite{Stein1971Weiss}, then polyconvolution $\underset{F_c , \L}{*}(f,g,h)$ belongs to $L_s (\R_+)$. To prove the inequality \eqref{eq5.10},  we choose the function $$k(x)={\rm sign}\big[(\underset{F_c,\L}{*}(f,g,h))(x)\big]^{s}.\big[(\underset{F_c,\L}{*}(f,g,h))(x)\big]^{s \textfractionsolidus s_1}.$$ Then $k \in L_{s_1} (\R_+)$, with the norm $\lv k \lv_{L_{ s_1} (\R_+)} = \lv \underset{F_c,\L }{*}(f,g,h)\lv^{s \textfractionsolidus s_1}_{L_{ s} (\R_+)}$. Applying inequality \eqref{eq5.1} to such  $k(x)$, we get
	
	$$
	\begin{aligned}\lv\underset{F_c,\L}{*}(f,g,h) \lv^{s}_{L_{s}(\R_+)}&= 
	\int\limits_{0}^{\infty}|(\underset{F_c,\L }{*}(f,g,h))(x)|^{s} dx = \bigg|\int\limits_{0}^{\infty}(\underset{F_c,\L}{*}(f,g,h))(x). \,k(x) dx\bigg|\\&\leq w^{(1-q)\textfractionsolidus q} v^{(1-r)\textfractionsolidus r} \lv f\lv_{L_p (\R_+)}\lv g\lv_{L_q (\R_+, (w+x)^{q-1})} \lv h\lv_{L_r (\R_+, (v+x)^{r-1})}\lv k \lv_{L_s (\R_+)}\\& = w^{(1-q)\textfractionsolidus q} v^{(1-r)\textfractionsolidus r} \lv f\lv_{L_p (\R_+)}\lv g\lv_{L_q (\R_+, (w+x)^{q-1})} \lv h\lv_{L_r (\R_+, (v+x)^{r-1})}\lv \underset{F_c,\L }{*}(f,g,h)\lv^{s \textfractionsolidus s_1}_{L_{ s} (\R_+)}.
	\end{aligned}$$
	\noindent Since $s - \frac{s}{s_1} =1$, from the above equality, we arrive at the conclusion of the corollary.
\end{proof}

\subsection{Norm estimation in the $L_p (\R_+ , \rho_j)$ weighted spaces}\label{subsection42}

\begin{theorem}\label{theorem6.1}
	Assume that $\rho_1, \rho_2, \rho_3$  are non-vanishing positive functions such that polyconvolution $\underset{F_c, \L}{*}(\rho_1, \rho_2,\rho_3)$ is well-defined. For any $F_j \in L_p (\R_+, \rho_j)$ with $p >1$,  the following $L_p (\R_+)$ weighted inequality holds true 
	\begin{equation}\label{eq6.1}
	\big\lv  \underset{F_c, \L}{*}(F_1 \rho_1, F_2 \rho_2, F_3 \rho_3)\underset{F_c , \L}{*}(\rho_1, \rho_2, \rho_3)^{\frac{1}{p}-1}\big\lv_{L_p (\R_+)}\leq \prod\limits_{j=1}^3\lv F_j\lv_{L_p (\R_+,\rho_j)},
	\end{equation}
	where $\underset{F_c ,L}{*}(.,.,.)$ is defined by formula \eqref{eq3.1} and $L_p (\R_+ , \rho_j)$ are weighted spaces with respect to a positive measure $\rho_j (x)dx$ equipped with the norm $\lv F_j \lv_{L_p (\R_+ , \rho_j)}=\bigg( \int\limits_0^{\infty} |F_j(x)|^p \rho_j (x) dx \bigg)^{1/p}<+ \infty, j=1,2,3.$
\end{theorem}
\begin{proof}
	Based on Definition \ref{dachapFcLaplace} of the polyconvolution for Fourier cosine, Laplace integral transforms, we get
	\begin{equation}\label{eq6.2}
	\begin{aligned}
	&\big\lv  \underset{F_c, \L}{*}(F_1 \rho_1, F_2 \rho_2, F_3 \rho_3)\underset{F_c , \L}{*}(\rho_1, \rho_2, \rho_3)^{\frac{1}{p}-1}\big\lv^p_{L_p (\R_+)}= \int\limits_{0}^{\infty} \big| (\underset{F_c ,\L}{*}(F_1 \rho_1, F_2 \rho_2, F_3 \rho_3))(x)\big|^p. \big|(\underset{F_c , \L}{*}(\rho_1, \rho_2, \rho_3))(x)\big|^{1-p} dx\\
	&=\frac{1}{\pi} \int\limits_{0}^{\infty}\biggl\{ \bigg| \int\limits_{R^3_+} \Phi(x,u,v,w).(F_1 \rho_1)(u)(F_2 \rho_2)(v)(F_3 \rho_3)(w)dudvdw\bigg|^p \times
	\bigg| \int\limits_{R^3_+} \Phi(x,u,v,w). \rho_1 (u)\rho_2 (v) \rho_3 (w)dudvdw\bigg|^{1-p}\biggr\}dx,
	\end{aligned}
	\end{equation}
	where $\Phi(x,u,v,w)$ is defined by \eqref{eq3.2}.
	On the other hand, using  H\"older's inequality for $q$ as the exponential conjugate to $p$, we obtain
	
	\begin{equation}\label{eq6.3}
	\begin{aligned}
	&\bigg| \int\limits_{R^3_+} \Phi(x,u,v,w)(F_1 \rho_1)(u)(F_2 \rho_2)(v)(F_3 \rho_3)(w)dudvdw\bigg|\\
	&\leq \biggl\{\int\limits_{R^3_+}|\Phi(x,u,v,w)|.|F_1 (u) |^p\ \rho_1 (u)|F_2 (v)|^p\ \rho_2 (v)|F_3 (w) |^p\ \rho_3 (w)dudvdw\biggr\}^{\frac{1}{p}}\times \biggl\{\int\limits_{R^3_+}|\Phi(x,u,v,w)|.\rho_1 (u)\rho_2 (v)\rho_3 (w)dudvdw \biggr\}^{\frac{1}{q}}
	\end{aligned}
	\end{equation}
	Deducing directly from the formulas \eqref{eq6.2} and \eqref{eq6.3}, we have
	\begin{equation*}
	\begin{aligned}
	&\big\lv  \underset{F_c, \L}{*}(F_1 \rho_1, F_2 \rho_2, F_3 \rho_3)\underset{F_c , \L}{*}(\rho_1, \rho_2, \rho_3)^{\frac{1}{p}-1}\big\lv^p_{L_p (\R_+)}\\&\leq\frac{1}{\pi} \int\limits_0^{\infty} \biggl\{ \bigg( \int\limits_{R^3_+} |\Phi(x,u,v,w)|. |F_1 (u)|^p\ \rho_1 (u) |F_2 (v) |^p\ \rho_2 (v)|F_3 (w) |^p\ \rho_3 (w)dudvdw\bigg)\times \\&\times \bigg( \int\limits_{R^3_+} |\Phi(x,u,v,w)|.\rho_1 (u)\rho_2(v)\rho_3(w)dudvdw\bigg)^{\frac{p}{q}}\times \bigg( \int\limits_{R^3_+} |\Phi(x,u,v,w)|.\rho_1 (u)\rho_2(v)\rho_3(w)dudvdw\bigg)^{1-p}\biggr\}dx.
	\end{aligned} 
	\end{equation*}
	Since $1\textfractionsolidus p+1\textfractionsolidus q=1$, then $\frac{p}{q} +1-p=0$. Therefore we infer that
	$\big\lv  \underset{F_c, \L}{*}(F_1 \rho_1, F_2 \rho_2, F_3 \rho_3)\underset{F_c , \L}{*}(\rho_1, \rho_2, \rho_3)^{\frac{1}{p}-1}\big\lv^p_{L_p (\R_+)}$ $\leq \frac{1}{\pi}\int\limits_{R^4_+} |\Phi(x,u,v,w)|. |F_1 (u)|^p\ \rho_1 (u) |F_2 (v) |^p\ \rho_2 (v)|F_3 (w) |^p\ \rho_3 (w) dudvdwdx.$
	By the assumption  $F_j \in L_p (\R_+ , \rho_j)$, using Fubini's theorem for the right-hand side of the above equality, we obtain
	$$\begin{aligned}
	&\big\lv  \underset{F_c, \L}{*}(F_1 \rho_1, F_2 \rho_2, F_3 \rho_3)\underset{F_c , \L}{*}(\rho_1, \rho_2, \rho_3)^{\frac{1}{p}-1}\big\lv^p_{L_p (\R_+)}\\
	&\leq \frac{1}{\pi} \bigg(\int\limits_{0}^{\infty} | \Phi(x,u,v,w)|dx \bigg)\bigg(\int\limits_{0}^{\infty} |F_1 (u) |^p\ \rho_1 (u)du\bigg)\bigg(\int\limits_{0}^{\infty} |F_2 (v) |^p\ \rho_2 (v)dv\bigg) \bigg(\int\limits_{0}^{\infty} |F_3 (w)|^p\ \rho_3 (w)dw\bigg)\\
	&=\frac{1}{\pi}.\pi \lv F_1\lv_{L_p (\R_+,\rho_1)}\lv F_2\lv_{L_p (\R_+,\rho_2)}\lv F_3\lv_{L_p (\R_+,\rho_3)}.
	\end{aligned}$$
	The theorem is proved.
\end{proof}
\noindent Taking $\rho_2(x)\equiv 1, \forall x\in\R_+$, and combining with Theorem \ref{theorem6.1} and \eqref{eq3.4}, we arrive at the following corollary.

\begin{corollary}\label{hqSaitoh}
	Let $\rho_1 \equiv 1, \forall x \in \R_+$  and 
	$0 \leq \rho_2 , \rho_3 \in L_1 (\R_+)$ such that $\underset{F_c ,\L}{*}(1,\rho_2 ,\rho_3)$ is well-defined. Then, for any functions $F_1 \in L_p (\R_+), F_2 \in L_p (\R_+, \rho_2), F_3 \in L_p (\R_+,\rho_3)$ with $ p>1$, we have the following estimate
	\begin{equation*}\label{eq6.5}
	\lv \underset{F_c ,\L}{*}(F_1,F_2\rho_2 ,F_3\rho_3)\lv_{L_p (\R_+)} \leq \lv \rho_2 \lv^{1-\frac{1}{p}}_{L_1 (\R_+)} \lv \rho_3 \lv^{1-\frac{1}{p}}_{L_1 (\R_+)} \lv F_1\lv_{L_p(\R_+)}\lv F_2\lv_{L_p(\R_+,\rho_2)}\lv F_3\lv_{L_p(\R_+,\rho_3)}.
	\end{equation*}
\end{corollary}
\noindent Here we consider an illustrative example for evaluation \eqref{eq6.5} as follows. Choosing $\rho_1(u)\equiv 1 ,\forall u \in \R_+$ and $\rho_2 (v)=e^{-v}, \rho_3 (w)=e^{-2w} \in L_1 (\R_+)$. Then using the formulas \eqref{eq3.1} and \eqref{eq3.5}, together with Fubini's theorem, we obtain
$$\begin{aligned}
&|(\underset{F_c ,\L}{*}(1,e^{-v},e^{-2w}))(x) |\leq \frac{1}{\pi } \int\limits_{R^3_+}|\Phi(x,u,v,w). e^{-v}e^{-2w}| dudvdw\\&=\frac{1}{\pi}\bigg(\int\limits_{0}^{\infty}| \Phi(x,u,v,w)|\ du\bigg)\bigg(\int\limits_{0}^{\infty}e^{-v}dv\bigg)\bigg(\int\limits_{0}^{\infty}e^{- 2w}dw\bigg)
= \lv e^{-v}\lv_{L_1 (\R_+)}\lv e^{-2w}\lv_ {L_1 (\R_+)}=\frac{1}{2}.
\end{aligned}$$
Thus polyconvolution $\underset{F_c ,\L}{*}(1,e^{-v},e^{-2w})$ is well-defined for all functions $F_1 \in L_p (\R_+) , F_2 \in L_p (\R_+ , e^{-v}),$ and $F_3 \in L_p (\R_+ , e^{-2w})$ with $p>1$. We obtain the following estimate
$$\lv \underset{F_c ,\L}{*}(F_1, F_2 e^{-v},F_3 e^{-2w})\lv_{L_p (\R_+)} \leq \frac{1}{2^{1-1\textfractionsolidus p}}\lv F_1\lv_{L_p (\R_+)}\lv F_2\lv_{L_p (\ R_+,e^{-v})}\lv F_3\lv_{L_p (\R_+,e^{-2w})}.$$
\section{Some applications}\label{apply}
The aim of this section is to consider the integral equation of polyconvolution type with
the Toeplitz plus Hankel kernels firstly posed in
\cite{Tsitsiklis1981Levy} and integro-differential equation of Barbashin type in \cite{Appell2000KalZabre}. By constructing $(F_c, \L)$-polyconvolutions, we obtain a necessary and sufficient condition for the solvability and unique explicit $L_p$-solutions $(p=1,2)$ respectively of those equations.
\subsection{On the Toeplitz plus Hankel type integral equation}
Following \cite{Tsitsiklis1981Levy}, we consider the integral equation of the form \begin{equation}\label{eq7.1}
f(x)+\int\limits_{0}^{T} \big[k_1 (x+y)+k_2 (x-y) \big]f(y)dy=\vp (x), \ 0\leq x,y \leq T,
\end{equation}where $\vp$ is a given function, $f$ an  unknown function and $K(x,y):= k_1 (x+y)+k_2 (x-y)$ is the kernel of the equation. Eq. \eqref{eq7.1}  with a Hankel
$k_1( x + y )$ or Toeplitz $k_2 ( x - y )$ kernel has attracted attention of many authors as they have practical applications in such diverse
fields such as scattering theory, fluid dynamics, linear filtering theory, inverse scattering problems in quantum mechanics, problems in radiative wave transmission, and further applications in medicine and biology (refer \cite{Agranovich1963Marchenko,Chadan1977Sabatier}). Eq. \eqref{eq7.1} has been carefully studied when $K(x,y)$ is a Toeplitz $k_2 (x - y)$ or Hankel kernel $k_1 (x + y)$. One notable case is when $K (x,y) = k_2 (x - y)+ k_2 ( x + y )$ , i.e. the Toeplitz and Hankel kernels  generated by the same function $k_2$ have been investigated in \cite{kagiwada1974integral}, namely by setting up Eq. \eqref{eq7.1}  in the domain $I = [0;T],$ $T > 0$ with kernel $K(x, \tau)=\int\limits_{0}^{1}\big[e^{-\frac{|x-\tau|}{\theta}}+e^{-\frac{x+\tau}{\theta}} r(\theta)\big] w(\theta)\ d \theta.$ In \cite{Tsitsiklis1981Levy}, Tsitsiklis and Levy considered Eq. \eqref{eq7.1} with general Toeplitz plus Hankel kernels $k_1 (x+y)+k_2 (x-y)$. This approach leads to Eq. \eqref{eq7.1} being a generalization
of Levinson's equation considered by Chanda and Sabatier in \cite{Chadan1977Sabatier} for kernel case is Toeplitz function based on the Gelfand-Levitan's method \cite{Gelfand1951Levitan} and the approach by Marchenko \cite{Agranovich1963Marchenko} for kernel is the Hankel function. 
A special case of Eq. \eqref{eq7.1} is the more extensive result of Tsitsiklis-Levy involving the Kontorovich-Lebedev method with two Toeplitz plus Hankel's type kernels defined in the infinite domain $I\equiv \R_+$ as follows: $
\lambda \int_{\mathbb{R}_{+}} K_1(x, \tau) f(\tau) d \tau+\mu \int_{\mathbb{R}_{+}} K_2(x, \tau) f(\tau) d \tau=g(x), x \in \mathbb{R}_{+},
$
where  constants $\lambda, \mu \in \mathbb{C}$ are predetermined, with the kernels $K_1(x, \tau), K_2(x, \tau)$  has form
$
K_1(x, \tau)=\frac{1}{\pi^2} \int_{\R_+}\left[\sinh (\tau-\theta) e^{-x \cosh (\tau-\theta)}+\sinh (\tau+\theta) e^{-x \cosh (\tau+\theta)}\right] \varphi(\theta) d \theta,$ and
$K_2(x, \tau)=\frac{1}{2 \pi x} \int_{\R_+}\left[e^{-x \cosh (\tau-\theta)}+e^{-x \cosh (\tau+\theta)}\right] h(\theta) d \theta,
$
has been proven in \cite{tuan2018hoanghong}.
However, the solution of \eqref{eq7.1} in closed-form for general case is still open.
Being different from other approaches, our idea is to reduce the original integral equation to become the linear equation by using the $(F_c,\L)$-polyconvolution. Thus, we will obtain the $L_1$-solution with the simultaneous help of the factorization properties and Winner-Lévy's Theorem \cite{NaimarkMA1972}. Namely,  for kernels $k_1 (x+y) = \frac{1}{\sqrt{2\pi}}g(x+y);\ k_2 (x-y) = \frac{1}{\sqrt{2\pi}}g(|x-y|)$ and choosing $\vp = \underset{F_c , \L}{*}(g,h,\xi)$, considering on infinite range $(0,T)\equiv (0,\infty)$, then the integral equation with Toeplitz plus Hankel kernel \eqref{eq7.1} can be rewritten in the convolution form 
\begin{equation}\label{eq7.2}
f(x)+(f\underset{F_c}{*}g)(x)= (\underset{F_c , \L}{*}(g,h,\xi))(x).
\end{equation}
\begin{theorem}\label{theorem7.1}
	Suppose that $h$ and $\xi$ are given functions belonging to  $L_1 (\R_+)$ and $g$ be $L_1$-Lebesgue integrable function over $\R_+$  satisfy the condition  $1+(F_c g)(y)\neq 0$ for any $y \in \R_+$. Then, equation \eqref{eq7.2}  has the unique solution in $L_1 (\R_+)$ which can be represented in the form $f(x)=(\underset{F_c , \L}{*}(\e,h,\xi))(x)$. Here $\e \in L_1 (\R_+)$ is defined by $(F_c \e)(y)=\frac{(F_c g)(y)}{1+(F_c g)(y)}$. Furthermore, the following $L_1$-norm estimate holds
	$\lv f\lv_{L_1 (\R_+)} \leq \lv \e\lv_{L_1 (\R_+)} \lv h\lv_{L_1 (\R_+)}\lv \xi\lv_{L_1 (\R_+)}.$
\end{theorem}	
\begin{proof}
	Applying the Fourier cosine transform to both sides of \eqref{eq7.2}, we get $(F_c f)(y)+F_c (f\underset{F_c}{*}g)(y)=F_c (\underset{F_c , \L}{*}(g,h,\xi))(y).$ Using the factorization properties \eqref{eq2.2},\eqref{eq3.4}, together with the condition $1+(F_c g)(y)\neq 0$, $\forall y\in \R_+$, we obtain $$(F_c f)(y)\big[ 1+(F_c g)(y)\big]=(F_c g)(y)(\L h)(y)(\L \xi)(y).$$
	This is equivalent to $(F_c f)(y)=\frac{(F_c g)(y)}{1+(F_c g)(y)}(\L h)(y)(\L \xi)(y)$. The Winner-Lévy's Theorem \cite{NaimarkMA1972} for the Fourier transform says that if $k \in L_1 (\R)$, then $1+(F k)(y) \neq 0$ for any $y \in \R$ is a necessary and sufficient condition for the existence of a function $\psi $  belonging to $L_1 (\R)$ such that $(F \psi)(y)=\frac{(F k)(y)}{1+(F k)(y)}$. This theorem still holds true for Fourier cosin transform  in the $L_1 (\R_+)$ \cite{Titchmarsh1986}.  Consequently, by the Winner-Lévy's Theorem, then $1+(F_c g)(y)\neq 0$ is a necessary and sufficient condition for the existence
	of a function $\e \in L_1 (\R_+)$ such that $(F_c \e)(y)=\frac{(F_c g)(y)}{1+(F_c g)(y)}$. This means that  $(F_c f)(y)=(F_c \e)(y)(\L h)(y)(\L \xi)(y)= F_c (\underset{F_c, \L}{*}(\e, h,\xi))(y).$
	Therefore  $f(x)= (\underset{F_c, \L}{*}(\e, h,\xi))(x)$ almost everywhere for any  $x\in \R_+$. Moreover, since $\e, h, \xi$ are functions belonging to the $L_1 (\R_+)$, by Theorem \ref{theorem31}, we deduce that  $f(x) \in L_1 (\R_+)$. Applying the inequality \eqref{eq3.3}, we obtain norm estimation of the solution on $L_1$ space as follows $\lv f \lv_{L_1 (\R_+)} \leq \lv \e \lv_{L_1 (\R_+)} \lv h \lv_{L_1 (\R_+)}\lv \xi \lv_{L_1 (\R_+)}$.		
\end{proof}

\begin{remark}\label{remark71}
	Let  $p,q,r,s >1$ be real numbers such that $1\textfractionsolidus p+1\textfractionsolidus q+1\textfractionsolidus r=2+1\textfractionsolidus s$. Applying the evaluation \eqref{eq5.10},  for any functions $f\in L_s (\R_+), \e\in L_p (\R_+), h \in L_q (\R_+,(w+) x)^{q-1})$ and $\xi \in L_r (\R_+, (v+x)^{r-1})$, we get the following solution estimate
	$$\lv f\lv_{L_s (\R_+)} \leq w^{(1-q)\textfractionsolidus q} v^{(1-r)\textfractionsolidus r}  \lv \e\lv_{L_p (\R_+)} \lv h\lv_{L_q (\R_+,(w+x)^{q-1}} \lv \xi \lv_{L_r (\R_+,(v+x)^{r-1}}.$$ Furthermore , without loss of generality, it can be assumed that $h =h_1 \a,$ and $\xi = \xi_1 \b$ with $\a, \b \in L_1 (\R_+)$ such that polyconvolution $\underset{F_c , \L}{*}(1,\a,\b)$ is well-defined. By using Corollary \ref{hqSaitoh}, $\forall f \in L_p (\R_+)$, $\e\in L_1 (\R_+) \cap L_p (\R_+), h_1 \in L_1 (\R_+,\a) \cap L_p (\R_+,\a),$ and $\forall \xi_1 \in L_1 (\R_+,\b) \cap L_p (\R_+,\b)$ with $p>1$, we get an estimate of the boundedness in weighted $L_p$ spaces for solution of the equation \eqref{eq7.2} as follows
	$$\lv f\lv_{L_p (\R_+)} \leq \lv \a \lv^{1-\frac{1}{p}}_{L_1 (\R_+ )} \lv \b \lv^{1-\frac{1}{p}}_{L_1 (\R_+)} \lv \e\lv_{L_p(\R_+ )}\lv h_1\lv_{L_p(\R_+,\a)}\lv \xi_1\lv_{L_p(\R_+,\b)}.$$	
\end{remark}
\noindent\underline{Comments:}    It can be seen that the complexity in this study is due to the kernel $K(x,y)$ generated by two discriminant functions $k_1$ (Hankel kernel) and $k_2$ (Toeplitz kernel) are defined on an infinite interval. Thus our approach to the analyzing of equation \eqref{eq7.1} derives from Fourier cosine and Laplace integral transforms (abbreviated as $F_c, \L$)  over infinite intervals combining the structure of convolutions appropriately constructed. We show a unique explicit solution on $L_1 (\R_+)$ via Winner-Lévy's theorem. These are in stark contrast to recent results of Castro et al in \cite{castro2019,castro2020new} when the kernel $K(x, y)$ is generated by two discriminant functions $p$ (Hankel kernel) and $q$ (Wiener-Hopf kernel) are $2\pi$-periodic defined on one finite interval and the solution is obtained on the Hilbert space $L_2([0, 2\pi])$ based on Shannon's sampling method.

To illustrate for Theorem \ref{theorem7.1} and Remark \ref{remark71}, we consider the following example.
\begin{example}
	Let $g(x)=\sqrt{\frac{\pi}{2}}e^{-x}, h(x)=\frac{1}{2}\sqrt{\frac{\pi}{2}}e^{-2x}$ and  $\xi(x)=\frac{1}{3}\sqrt{\frac{\pi}{2}}e^{-3x}$. It is easy to check that $g,h,\xi$ are functions belonging to $ L_1 (\R_+)$. According to \cite{bateman1954}, formula 1.4.1, page 14, we have $F_c (\sqrt{\frac{\pi}{2}}e^{-x})=\frac{1}{1+y^2},$ then  $1+(F_c g)(y)=1+\frac{1}{1+y^2} \neq 0.$ This implies that $g(x)$ actually satisfies the condition of Theorem \ref{theorem7.1}. 
	Furthermore, we have $(F_c \e)(y)=\frac{(F_c g)(y)}{1+(F_c g)(y)}=\frac{1}{2+y^2},$ then $\e(x)=\frac{\sqrt{\pi}}{2}e^{-\sqrt{2}x}$ belongs to  $L_1 (\R_+)$. Therefore, the solution Eq. \eqref{eq7.2} belongs to $L_1 (\R_+)$ and is represented in the following form $f(x)=\big(\underset{F_c , \L}{*}(\frac{\sqrt{\pi}}{2}e^{-\sqrt{2}t}, \frac{\sqrt{\pi}}{2\sqrt{2}}e^{-2t},\frac{\sqrt{\pi}}{3\sqrt{2}}e^{-3t} )\big)(x).$ And we obtain the norm estimate in $L_1$ of the solution as follows
	$$\lv f\lv_{L_1 (\R_+)}\leq \big\lv \frac{\sqrt{\pi}}{2}e^{-\sqrt{2}t}\big\lv_{L_1 (\R_+)}\big\lv \frac{\sqrt{\pi}}{2\sqrt{2}}e^{-2t}\big\lv_{L_1 (\R_+)}\big\lv \frac{\sqrt{\pi}}{3\sqrt{2}}e^{-3t}\big\lv_{L_1 (\R_+)}=\frac{\pi\sqrt{\pi}}{134\sqrt{2}}.$$
	On the other hand, $h(x)=\frac{1}{2}\sqrt{\frac{\pi}{2}}e^{-2x}=\frac{1}{2}\sqrt{\frac{\pi}{2}}e^{-x}.e^{-x}$ and $\xi(x)=\frac{1}{3}\sqrt{\frac{\pi}{2}}e^{-3x}=\frac{1}{3}\sqrt{\frac{\pi}{2}}e^{-2x}.e^{-x}$. So we choose $h_1 =\frac{1}{2}\sqrt{\frac{\pi}{2}}e^{-x}$ and $ \xi_1=\frac{1}{3}\sqrt{\frac{\pi}{2}}e^{-2x}$ with $\a=\b=e^{-x} \in L_1 (\R_+)$. Using formulas \eqref{eq3.1}, \eqref{eq3.5} we obtain 
	
	$$\begin{aligned}
	(\underset{F_c , \L}{*}(1,\a,\b))(x)&=(\underset{F_c , \L}{*}(1,e^{-t},e^{-t}))(x):=\frac{1}{\pi}\int\limits_{R^3_+}\Phi(x,u,v,w)e^{-v}e^{-w}dudvdw\\
	&=\frac{1}{\pi}\bigg(\int\limits_{0}^{\infty} \Phi(x,u,v,w)du\bigg)\bigg(\int\limits_{0}^{\infty} e^{-v}dv\bigg)\bigg(\int\limits_{0}^{\infty} e^{-w}dw\bigg)=\frac{1}{\pi}. \pi\lv e^{-v}\lv_{L_1 (\R_+)}\lv e^{-w}\lv_{L_1 (\R_+)}=1.
	\end{aligned}$$
	Hence polyconvolution $\underset{F_c , \L}{*}(1,e^{-t},e^{-t})$ is well-defined with the  given functions $\a , \b$. And we have an estimate of $L_p$-solution for Eq. \eqref{eq7.2} based on Remark \ref{remark71} is 
	$$\lv f\lv_{L_p (\R_+)}\leq \big( \lv e^{-x}\lv^{1-\frac{1}{p}}_{L_1 (\R_+)}\big)^2 \big\lv \frac{\sqrt{\pi}}{2}e^{-\sqrt{2}x}\big\lv_{L_p (\R_+, e^{-x})}\big\lv \frac{\sqrt{\pi}}{3\sqrt{2}}e^{-2x}\big\lv_{L_p (\R_+,e^{-x})} = \frac{\pi\sqrt{\pi}}{24}\biggl\{\frac{1}{\sqrt{2}p(p+1)(2p+1)}\biggr\}^{\frac{1}{p}}.
	$$
\end{example}
\subsection{On linear integro-differential equations of Barbashin type}
The subsection is
concerned with integro-differential equation of the form 
\begin{equation}\label{eqBarbashin}
\frac{\p f(t,s)}{\p t}=c(t,s)f(t,s)+\int\limits_{a}^b K(t,s,\rho) f(t,s)d\rho +g(t,s).
\end{equation}
Here $c: J \times[a, b] \longrightarrow \mathbb{R}, K: J \times[a, b] \times[a, b] \longrightarrow \mathbb{R}$, and mapping $g: J \times[a, b] \longrightarrow \mathbb{R}$ are given functions, where $J$ is a bounded or unbounded interval, the function $f$ is unknown. The equation \eqref{eqBarbashin} was first studied  by E.A. Barbashin \cite{Barbashin1957} and his pupils. For this reason, this is nowadays called integro-differential  equation of Barbashin type or simply the Barbashin equation. Eq. \eqref{eqBarbashin} has been applied to many fields such as mathematical physics, radiation propagation, mathematical biology and transport problems, e.g., more details refer \cite{Appell2000KalZabre}.  One
of the characteristics of Barbashin equation is that studying solvability of the equation is
heavily dependent on the kernel $K(t,s,\rho)$ of the equation. In many cases, we can reduce Eq. \eqref{eqBarbashin} to the form of an ordinary differential equation and use the Cauchy integral operator or evolution operator to study it when the kernel does not depend on  $t$.  A typical example is  the stationary integro-differential equation of Barbashin type $\partial f(t, s) \textfractionsolidus \partial t=c(s) f(t, s)+\int\limits_{a}^{b} K(s, \rho) f(t, \rho) d \rho +g(t, s)$ i.e. the kernel $K$ does not depend on parameter $t$. In some other cases, we need to use the partial integral operator to study this equation (see \cite{Appell2000KalZabre}).  However, in the general case of $K(t,s,\rho)$ as an arbitrary kernel, the problem of finding a solution for Barbashin equation remains open. On the other hand, if we view ${A}$ as the operator defined by ${A}:=\partial\textfractionsolidus \partial t -c(t,s) \mathcal{I}$, where $\mathcal{I}$ is the identity operator, then Eq. \eqref{eqBarbashin}  is written in the following form
\begin{equation}\label{eq7.6}
{A}f(t,s)=\int\limits_a^b K(t,s,\rho)f(t,s)d\rho +g(t,s).
\end{equation}
\noindent To give a necessary and sufficient condition for the solvability and unique explicit $L_2$-solutions of Eq. \eqref{eq7.6}, first we need some following auxiliary lemmas.
\begin{lemma}\label{lemma7.1}
	Let $\vp$ be a function belonging to $L_2 (\R_+)$ and satisfy  $(1+y^2)|(F_c \vp)(y)|$ bounded $ \forall y>0$. For all $ f \in L_2 (\R_+)$, we have the following assertion
	\begin{equation}\label{eq7.4}
	F_c \big[ D(f\underset{F_c}{*}\vp)(t)\big](y)=(1+y^2)(F_c f)(y)(F_c \vp)(y),
	\end{equation}
	where $D$ is the second order differential operator having the form $D=\left(I-\frac{d^2}{dt^2}\right)$ and $(\underset{F_c}{*})$ is determined by \eqref{tichchapFc}.
\end{lemma}
\begin{proof}
	The argument is similar to that in the proof of  Lemma\ref{bode4.1}, it means if there are $f(y), yf(y), y^2 f(y) \in L_2 (\R_+)$, then $(F_c f)(t), \frac{d}{dt}(F_c f)(t), \frac{d^2}{dt^2}(F_c f)(t) \in L_2 (\R_+)$. Moreover $\frac{d^2}{dt^2}(F_c f)(t)=F_c (-y^2 f(y))(t)$ belongs to $L_2 (\R_+)$, implying that $$\big(1-\frac{d^2}{dt^2}\big)\big[(F_c f)(y)\big](t)=F_c \big[(1+y^2)f(y)\big](t).$$ On the other hand, we know that for all functions $f, \vp \in L_2 (\R_+)$, then $(f\underset{F_c}{*}\vp)(x)=F_c \big[ (F_c f)(y)(F_c \vp)(y)\big]$ (refer \cite{MusallamVKt2000}). Since $f \in L_2 (\R)$ then $(F_c f)\in L_2 (\R)$, together with assuming $(1+y^2)|(F_c \vp)(y)|$ is bounded $ \forall y>0$, we deduce that $(1+y^2)|(F_c \vp)(y)|$ is a bounded function on $\R$. Therefore  $(1+y^2)(F_c \vp)(y)(F_c f)(y)\in L_2 (\R_+)$. From the foregoing, we obtain
	$$\big(1-\frac{d^2}{dt^2}\big)(f\underset{F_c}{*}\vp)(t)=\big(1-\frac{d^2}{dt^2}\big)F_c \big[(F_c f)(y)(F_c \vp)(y) \big](t)= F_c \big[(1+y^2)(F_c f)(y)(F_c \vp)(y)\big](t) \in L_2 (\R_+).$$
	Furthermore, $F_c : L_2 (\R_+)\longleftrightarrow L_2 (\R_+)$ is unitary transformation \cite{Titchmarsh1986}, so the  expression \eqref{eq7.4} is directly deduced.
\end{proof}

\noindent Now, by putting operator $Af(t,s):=D(f\underset{F_c}{*}\vp)(t)$ and choosing the kernel as $K(t,s,\rho)=\frac{-1}{\pi}\int\limits_{R^2_+} \Phi (t,s,\rho, w)$ $\eta(\rho)\xi(w)d\rho dw,$ $ \forall t>0$, and considering $(a,b)=(0,\infty)$, then Eq. \eqref{eq7.6} can be rewritten in the following form 
\begin{equation}\label{eq7.8}
D(f\underset{F_c}{*}\vp)(t) + (\underset{F_c , \L}{*}(f,\eta, \xi))(t)=g(t),\ t>0.\end{equation} Here $(f\underset{F_c}{*}\vp)$ and $\underset{F_c , \L}{*}(f,\eta, \xi)$ are determined by \eqref{tichchapFc}; \eqref{eq3.1} respectively, and $f$ is the unknown function needed to find.
\begin{theorem}\label{theorem7.2}
	Let $g$ belong to $L_2 (\R_+)$ and $\vp \in L_2 (\R_+)$ such that $(1+y^2)|(F_c \vp)(y)|$ is bounded $\forall y >0$. Suppose that $\eta, \xi \in L_1 (\R_+)$ are given functions and satisfy $(1+y^2)(F_c \vp)(y)+(\L \eta)(y)(\L \xi)(y) \neq 0$ for any $y \in \R_+$.
	Then equation \eqref{eq7.8}	has the unique solution in $L_2 (\R_+)$.
\end{theorem}
\begin{proof}
	Applying the Fourier cosine transform to both sides of \eqref{eq7.8}, we get $$F_c \big[ D(f\underset{F_c}{*}\vp)(t)\big](y)+F_c \big[ (\underset{F_c , \L}{*}(f,\eta, \xi))(t)\big](y)=(F_c g)(y),\ \forall y>0.$$
	Inferred from using the factorization property \eqref{eq3.4} and equality \eqref{eq7.4}, we obtain
	$$(F_c f)(y)\big[(1+y^2)(F_c \vp)(y)+(\L \eta)(y)(\L \xi)(y) \big]=(F_c g)(y).$$
	On the left-hand side of the above equality, since $\eta, \xi$ belong to the $L_1 (\R_+)$ then $(\L\eta )(y)$ and $(\L\xi )(y)$  are finite. Together with condition of theorem as $(1+y^2)|(F_c \vp)(y)|< \infty$, we deduce that $0 \neq (1+y^2)(F_c \vp)(y)+(\L \eta)(y)(\L \xi)(y)$ is a bounded function on $\R_+$. Moreover, for all $g \in L_2 (\R_+)$ then $(F_c g) \in L_2 (\R_+)$. Hence $$(F_c f)(y)=(F_c g)(y)\textfractionsolidus [(1+y^2)(F_c \vp)(y)+(\L \eta)(y)(\L \xi)(y)] \in L_2 (\R_+),$$ implying that solution $f$ of the Eq. \eqref{eq7.8} belongs to  $L_2 (\R_+)$. Based on the inverse Fourier cosine transform, we represent explicitly the solution in the following form 
	$f(t)=\sqrt{\frac{2}{\pi}}\int\limits_{0}^{\infty}\frac{(F_c g)(y)}{(1+y^2)(F_c \vp)(y)+(\L \eta)(y)(\L \xi)(y)}\cos(ty)dy\ $ almost everywhere on $\R_+$ (see \cite{Sogge1993fourier,WRudin1987}).\end{proof}

To illustrate for Theorem \ref{theorem7.2}, we consider the following example.
\begin{example}
	Let $g(x)=\vp(x) =\sqrt{\frac{\pi}{2}}e^{-x} \in L_2 (\R_+)$, then we have $F_c \big(\sqrt{\frac{\pi}{2}}e^{-x}\big)=\frac{1}{1+y^2}$, implying that $(1+y^2)|(F_c \vp)(y)|=1$ is bounded $\forall y >0$. Now, we choose $\eta =e^{-x}, \xi=xe^{-x} \in L_1 (\R_+)$. According to the formulas 4.5.1, page 143 and  4.5.2, page 144 in \cite{bateman1954}, we obtain $\L(e^{-x})=\frac{1}{1+y}$ and $\L(xe^{-x})=\frac{1}{(1+y)^2}$ finite. Therefore
	$$(F_c f)(y)=\frac{F_c \big(\sqrt{\frac{\pi}{2}}e^{-x}\big)}{(1+y^2)F_c \big(\sqrt{\frac{\pi}{2}}\big)e^{-x}+\L (e^{-x})\L(xe^{-x})}=\frac{(1+y)^3}{\big[(1+y)^3 +1\big](1+y^2)}\in L_2 (\R_+).$$
	This means that equation \eqref{eq7.8} has a unique solution in $L_2 (\R_+)$ and is represented as
	$$f(t)=\sqrt{\frac{2}{\pi}}\int\limits_{0}^{\infty}\frac{(1+y)^3}{(1+y^2)\big[(1+y)^3 +1\big]}\cos(ty)dy,$$ with $t>0.$ And we get an estimate $|f(t)|\leq \sqrt{\frac{2}{\pi}}\operatorname{Arc}\tan y\big|_{0}^{\infty}=\sqrt{\frac{\pi}{2}}$.	
\end{example}
For the case if we put operator $\mathcal{A} f(t,s):=D(\underset{F_c , \L}{*}(f,\eta,\xi))(t)$ and choose the kernel as $K(t,s,\rho)=\frac{-1}{\sqrt{2\pi}}\big[h(t+s)+h(|t-s|)\big].$ Then Eq. \eqref{eq7.6} becomes a type of Barbashin equation with Toeplitz plus Hankel kernels as follows
\begin{equation}\label{eq7.11}
D(\underset{F_c , \L}{*}(f,\eta,\xi))(t)+\frac{1}{\sqrt{2\pi}} \int\limits_{0}^\infty \big[h(t+s)+h(|t-s|)\big] f(s)d(s) = g(t),\ t>0, 
\end{equation}
\begin{theorem}\label{theorem7.3}
	Let $ g, h \in L_2 (\R_+)$, and $\eta, \xi$ are  $L_1$-Lebesgue integrable function over $\R_+$  such that $(1+y^2)\ |(\L \eta)(y)(\L \xi)(y)|$ is bounded $\forall y>0$, and satisfies the codition
	\begin{align}\label{dieukienTheorem7.3}
	(F_c g)(y)\textfractionsolidus \big[(1+y^2)(\L \eta)(y)(\L \xi)(y)+(F_c h)(y)\big] \in L_2 (\R_+). 
	\end{align} Then equation \eqref{eq7.11} 
	has the unique solution in $L_2 (\R_+)$. 
\end{theorem}

We need the following Lemma. 
\begin{lemma}\label{lemma7.2}
	If $\eta , \xi \in L_1 (\R_+)$ such that $(1+y^2)|(\L \eta)(y)(\L \xi)(y)|$ is a bounded function. Then for all $f \in L_2 (\R_+)$, we have
	\begin{equation}\label{eq7.5}
	F_c \big[D(\underset{F_c , \L}{*}(f,\eta,\xi))(t) \big](y)=(1+y^2)(F_c f)(y)(\L \eta)(y)(\L \xi)(y),\ \forall y>0,
	\end{equation} 
	where $D$ is the second order differential operator having the form $D=\big(I-\frac{d^2}{dt^2}\big)$ and $(F_c,\L)$-polyconvolution is determined by \eqref{dachapFcLaplace}.
	The proof of this Lemma is similar to that done with Lemma \ref{lemma7.1}.
\end{lemma} 
\begin{proof}\textit{of Theorem \ref{theorem7.3}.}\\
	Following the definition of generalized convolution \eqref{tichchapFc}, we convert \eqref{eq7.11} becoming
	$D(\underset{F_c , \L}{*}(f,\eta,\xi))(t)+(f \underset{F_c}{*}h)(t)=g(t),$ $t>0$.
	Applying the Fourier cosine transform to both sides of this equation, we have 
	$$F_c \big[D(\underset{F_c , \L}{*}(f,\eta,\xi))(t)\big](y)+F_c(f \underset{F_c}{*}h)(t)= (F_c g)(y),\ \forall y>0.$$
	Using \eqref{eq7.5} and factorization property \eqref{eq2.2}, we obtain $(F_c f)(y)\big[ (1+y^2)(\L \eta)(y)(\L\xi)(y)+(F_c h)(y)\big]=(F_c g)(y).$ This means that $(F_c g)(y)\textfractionsolidus \big[(1+y^2)(\L \eta)(y)(\L \xi)(y)+(F_c h)(y)\big] = (F_c f)(y)$  belongs to $L_2 (\R_+)$ due to condition \eqref{dieukienTheorem7.3}. Therefore, the solution $f$ of Eq. \eqref{eq7.11} belongs to $L_2 (\R_+)$. Through the inverse of Fourier cosine transform, we can rewrite explicitly the solution  as follows $$f(t)=\sqrt{\frac{2}{\pi}}\int\limits_{0}^{\infty}\frac{(F_c g)(y)}{(1+y^2)(\L \eta)(y)(\L \xi)(y)+(F_c h)(y)}\cos(ty)dy$$	almost everywhere on $\R_+$ (see \cite{WRudin1987,Sogge1993fourier}).
\end{proof}
To illustrate for Theorem \ref{theorem7.3}, we consider the following example.
\begin{example}
	We choose $\eta(x) = \xi(x) =e^{-x} \in L_1 (\R_+)$ and let $g(x) =h(x) =\sqrt{\frac{\pi}{2}}e^{-x} \in L_2 (\R_+)$. Hence, we have  $\L(e^{-x})(y)=\frac{1}{1+y}$ and $F_c \big(\sqrt{\frac{\pi}{2}}e^{-x}\big)=\frac{1}{1+y^2}$. Then $(1+y^2)|(\L e^{-x})(y)(\L e^{-x})(y) |=\frac{1+y^2}{(1+y)^2}\leq 1$, implying that $\eta,\xi$ satisfy the bounded condition. We have $(F_c f)(y)=\frac{F_c \big(\sqrt{\frac{\pi}{2}}e^{-x}\big)(y)}{(1+y^2)(\L e^{-x})(y)(\L e^{-x})(y)+F_c \big(\sqrt{\frac{\pi}{2}}e^{-x}\big)(y)}=\frac{(1+y)^2}{(1+y^2)^2 + (1+y)^2} \in L_2 (\R_+),$ this proves that the given functions actually satisfy condition \eqref{dieukienTheorem7.3}. Then the Eq. \eqref{eq7.11} has a unique solution in $L_2 (\R_+)$ and expressed by $f(t)=\sqrt{\frac{2}{\pi}}\int\limits_{0}^{\infty}\frac{(1+y)^2}{(1+y^2)^2 +(1+y)^2}\cos(ty)dy$ with $t>0,$ and we get the evaluation $|f(t)|\leq 2 \sqrt{\frac{2}{\pi}}\operatorname{Arc}\tan y\big|_{0}^{\infty}=\sqrt{2\pi}$.
\end{example}
\subsection{On differential equation}
The Theorem 4.1 in \cite{KilbasSrivastavaTrujillo2006}, Chapter 4, page 224, gives a closed-form solution of the Cauchy-type problem. 
\begin{align}\label{Cauchyproblem}\begin{cases}
&\left(D_{a+}^{\alpha} f\right)(x)-\lambda f(x)=\Phi(x), \quad(a<x \leq b;\ \alpha>0 ; \lambda \in \mathbb{R}) \\
&\left(D_{a+}^{\alpha-k} f\right)(a+)=b_{k}, \quad\left(b_{k} \in \mathbb{R} ; k=1, \cdots, n=-[-\alpha]\right),
\end{cases}\end{align}
where  $\Phi(x) \in C_{\gamma}[a, b],\ (0 \leq \gamma<1)$ with the Riemann-Liouville fractional derivative $\left(D_{a+}^{\alpha} f\right)(x)$ of order $\alpha >0$ are given by
$\left(D_{a+}^{\alpha} f\right)(x):=\frac{1}{\Gamma(n-\alpha)}\left(\frac{d^n}{d x^n}\right) \biggl\{\int\limits_{a}^{x} \frac{f(t)}{(x-t)^{\alpha-n+1}} d t\biggr\}, (n=[\alpha]+1 ; x>a).$
In this part, we consider the closed-form solutions of a narrow class of differential equations instead of the equation in the problem \eqref{Cauchyproblem} by choosing $\lambda = -1$ and  substituting operator $\left(D_{a+}^{\alpha} f\right)$ by   $D\big(\underset{F_c, \L}{*}(f, \eta,\xi)\big)$ where $D$ is the second-order differential operator. Then, the equation in  Cauchy-type problem \eqref{Cauchyproblem} can be rewritten 
\begin{equation}\label{ptviphan}
f(x) +(\T_{ \eta, \xi} f)(x) = g(x),\ x>0,
\end{equation}
where $\T_{ \eta, \xi}$ is defined by \eqref{eq4.1}.
\begin{theorem}\label{Theorem7.4}
	Let $g \in L_2 (\R_+)$. Suppose that $\eta, \xi$ are functions belonging to $\A (\R_+)$ (see Definition \ref{dnAR+}), such that the condition $0 \neq 1+(1+y^2)(\L \eta)(y)(\L \xi)(y)< \infty$. Then the Eq. \eqref{ptviphan} has a unique solution in $L_2 (\R_+)$ which can be presented in the form $f(x)=\sqrt{\frac{2}{\pi}}\int\limits_{0}^{\infty}\frac{(F_c g)(y)}{ 1+(1+y^2)(\L \eta)(y)(\L \xi)(y)} \cos(xy)dy,\ x>0.$ 
\end{theorem}
\begin{proof}
	Based on \eqref{eq4.1}, Eq. \eqref{ptviphan} is converted into the following form $f(x)+\big(1-\frac{d^2}{dx^2}\big)(\underset{F_c , \L}{*}(f,\eta, \xi))(x)=g(x)$, $ x>0.$ Applying the Fourier cosine $(F_c)$ transform to both sides of the above equation, we obtain
	$$(F_c f)(y)+F_c \big[ \big(1-\frac{d^2}{dx^2}\big)(\underset{F_c , \L}{*}(f,\eta, \xi))(x) \big](y)=(F_c g)(y), \ y>0.$$
	Using formula \eqref{eq4.7}, we obtain $(F_c f)(y)+(1+y^2)(F_c f)(y)(\L \eta)(y)(\L \xi)(y)=(F_c g)(y)$. This is equivalent to $(F_c f)(y)=\frac{(F_c g)(y)}{1+(1+y^2)(\L \eta)(y)(\L \xi)(y)}$ because the denominator of expression is non-zero under the condition of theorem. Since $g\in L_2 (\R_+)$ then $(F_c g )\in L_2 (\R_+)$. From the assumption $1+(1+y^2)(\L \eta)(y)(\L \xi)(y)$ is finite, we deduce that $\frac{1}{1+(1+y^2)(\L \eta)(y)(\L \xi)(y)}$ is a bounded function. This yields $(F_c g)(y).\frac{1}{1+(1+y^2)(\L \eta)(y)(\L \xi)(y)} \in L_2 (\R_+),$ implying that $(F_c f)(y) \in L_2 (\R_+)$ and $f$ belongs to $L_2 (\R_+)$ almost everywhere \cite{Sogge1993fourier}. According to the inverse formula of Fourier cosine transform, we get the solution in the explicit form as the conclusion of the theorem.
\end{proof}
We will end this article with an example illustrating the Theorem \ref{Theorem7.4}.
\begin{example}
	Let $\eta(x) =e^{ix}, \xi(x) =e^{-ix}$ and  $g(x)=\sqrt{\frac{2}{\pi}}K_0 (x)$, where $K_0 (y)=\int\limits_0^{\infty} \frac{1}{\sqrt{1+x^2}} \cos(xy)dx$ is a modified Bessel's function of the second kind. Following the formula 1.2.17, page 9 in \cite{bateman1954}, we obtain $(\L e^{ix})(y)=\frac{1}{y-i}$  and  $(\L e^{-ix})(y)=\frac{1}{y+i}$ are finite, hence $e^{\pm ix}\in \mathscr{A}(\R_+)$ and $1+(1+y^2)(\L e^{ix})(y)(\L e^{-ix})(y)=2$. On the other hand, $\int\limits_{0}^{\infty}|g(x)|^2 dx=\frac{\pi}{2}$, then $g(x)=\sqrt{\frac{2}{\pi}}K_0 (x)$ belongs to $L_2 (\R_+)$ and $(F_c g)=\frac{1}{\sqrt{1+y^2}}$. Thus, the functions given above completely satisfy the conditions of Theorem \ref{Theorem7.4} and we conclude that, Eq. \eqref{ptviphan} has a unique solution in $L_2 (\R_+)$ in the following form $f(x)=\frac{1}{\sqrt{2\pi}}K_0(x)\in L_2 (\R_+)$.\end{example}
\vskip 0.4cm
\noindent{\textbf{ACKNOWLEDGEMENTS}}\\
The author is grateful to Professor V\~u Kim Tu\^a\'n for useful discussions.  We would also like to thank all referees for carefully reading the manuscript and providing insightful comments that truly improved this paper.
\vskip 0.5cm
\noindent \textbf{CONFLICT OF INTEREST}\\
This work does not have any conflict of interest.
\vskip 0.4cm
\noindent \textbf{ORCID}\\
\noindent \textit{Trinh Tuan} {\color{blue} \url{https://orcid.org/0000-0002-0376-0238}}


\begin{thebibliography}{99}
\bibitem{delsarte1939} J. Delsarte, Une extension nouvelle de la théorie des fonctions presque-périodiques de Bohr, \textit{Acta Math.} \textbf{69} (1938), 259–317. \url{https://doi.org/10.1007/BF02547715}
\bibitem{BMLevitan1962} B. M. Levitan, Generalized translation operators and some of their applications, 1962. (in Russian). Israel Program for Scientific Translations, Nauka, Moscow.
\bibitem{Povzner1948} A. Povzner, On differential equations of Sturm-Liouville type on a half-axis. (in Russian), \textit{Matematicheskii Sbornik} \textbf{65} (1948), no. 1, 3–52.
\bibitem{kakichev1967convolution} V. A. Kakichev, On the convolution for integral transforms. (in Russian), \textit{Izv. Acad. Nauk BSSR. Ser. Fiz. Mat. Nauk.} \textbf{22} (1967), 48–57.
\bibitem{Kakichev1997}V. A. Kakichev, \textit{Polyconvolutions. Definitions, examples, convolutional equations,} Lecture notes. (in Russian), Taganrog State University of Radio Engineering, Taganrog, 1997.
\bibitem{Britvina2004}L. E. Britvina, On polyconvolutions generated by the Hankel transform, \textit{Math. Notes.} \textbf{76} (2004), no. 1, 18–24. \url{https://doi.org/10.1023/B:MATN.0000036738.90655.d5} 
\bibitem{britvina2007}L. E. Britvina, Generalized convolutions for the Hankel transform and related integral operators, \textit{Math. Nachr.} \textbf{280} (2007), no. 9-10, 962–970. \url{https://doi.org/10.1002/mana.200510528}
\bibitem{Thao2010Virchenko}N. X. Thao and N. A. Virchenko, On the polyconvolution for the Fourier cosine, Fourier sine, and Kontorovich–Lebedev integral transforms, \textit{Ukrainian. Math. J.} \textbf{62} (2011), no. 10, 1611–1624. \url{https://doi.org/10.1007/s11253-011-0453-8}
\bibitem{Sogge1993fourier}C. D. Sogge, \textit{Fourier integrals in classical analysis}, volume 105 of Cambridge Tracts in Mathematics, Cambridge University Press, 1993. \url{https://doi.org/10.1017/CBO9780511530029}
\bibitem{bateman1954}H. Bateman and A. Erdélyi, \textit{Tables of integral transforms}, McGraw-Hill Book Company, Inc, New York (NY), London, 1954.
\bibitem{Sneddon1972}I. N. Sneddon, \textit{The use of integral transform}, McGraw-Hill Book Company, Inc, New York (NY), London, 1972.
\bibitem{WRudin1987}W. Rudin, \textit{Real and complex analysis}, 3rd ed., McGraw-Hill International Book Company, New York, 1987.
\bibitem{VKTtuan2012}V. K. Tuan and T. Tuan, A real-variable inverse formula for the Laplace transform, \textit{Integr. Trans. Spec. Funct.} \textbf{23} (2012), no. 8, 551–555. \url{https://doi.org/10.1080/10652469.2011.609817}
\bibitem{Churchill1941}R. V. Churchill, \textit{Operational mathematics}, McGraw-Hill Science, Engineering \& Mathematics, New York, 1972.
\bibitem{MusallamVKt2000}F. Al-Musallam and V. K. Tuan, Integral transforms related to a generalized convolution, \textit{Result. Math.} \textbf{38} (2000), no. 3, 197–208. \url{https://doi.org/10.1007/BF03322007}
\bibitem{Titchmarsh1986}E. C. Titchmarsh, \textit{Introduction to the theory of Fourier integrals}, 3rd ed., Chelsea Publishing Co, New York, 1986.
\bibitem{Debnath2006Bhatta}L. Debnath and B. Bhatta, \textit{Integral transforms and their applications}, 2nd ed., Chapman \& Hall/CRC Press, New York, 2006. \url{https://doi.org/10.1201/9781420010916}
\bibitem{WatsonGN1933}G. N. Watson, General transforms, \textit{Proc. London Math. Soc.} \textbf{35} (1933), no. 2, 156–199. \url{https://doi.org/10.1112/plms/s2-35.1.156}
\bibitem{VKT1999JMAA}V. K. Tuan, Integral transforms of Fourier cosine convolution type, \textit{J. Math. Anal. Appl.} \textbf{229} (1999), no. 2, 519–529. \url{https://doi.org/10.1006/jmaa.1998.6177}
\bibitem{tuan2022Mediterranean}T. Tuan, Operational properties of the Hartley convolution and its applications, \textit{Mediterr. J. Math.} \textbf{19} (2022), 266. \url{https://doi.org/10.1007/s00009-022-02173-5}
\bibitem{Tuan2022MMA}T. Tuan, Some results of Watson and Plancherel–type integral transforms related to the Hartley, Fourier convolutions and applications, \textit{Math. Methods Appl. Sci.} \textbf{45} (2022), no. 17, 11158–11180. \url{https://doi.org/10.1002/mma.8442}
\bibitem{tuan2020Ukrainian}T. Tuan, On the Fourier-sine and Kontorovich–Lebedev generalized convolution transforms and their applications, \textit{Ukrainian Math. J.} \textbf{72} (2020), no. 2, 302–316. \url{https://doi.org/10.1007/s11253-020-01782-1}
\bibitem{YoungWH1912}W. H. Young, On classes of summable functions and their Fourier series, \textit{Proc. R. Soc. London, Ser. A.} \textbf{87} (1912), 225–229. \url{https://doi.org/10.1098/rspa.1912.0076}
\bibitem{AdamsFournier2003sobolev}R. A. Adams and J. F. Fournier, \textit{Sobolev spaces}, 2nd ed., Pure and Applied Mathematics, Vol. \textbf{140}, Elsevier/Academic Press, Amsterdam, 2003.
\bibitem{Saitoh2000}S. Saitoh, Weighted 
$L_p$-norm inequalities in convolutions, \textit{Survey on Classical Inequalities}, T. M. Rassias, (ed.), Mathematics and its Applications, Vol.  517, Springer, Dordrecht, 2000, pp. 225–234. \url{https://doi.org/10.1007/978-94-011-4339-4_8}
\bibitem{hoangtuan2017thaovkt}P. V. Hoang, T. Tuan, N. X. Thao, and V. K. Tuan, Boundedness in weighted  $L_p$
spaces for the Kontorovich–Lebedev–Fourier generalized convolutions and applications, \textit{Integr. Trans. Spec. Funct.} \textbf{28} (2017), no. 8, 590–604. \url{https://doi.org/10.1080/10652469.2017.1330825}
\bibitem{TuanVKTuan2023ITSF}T. Tuan and V. K. Tuan, Young inequalities for a Fourier cosine and sine polyconvolution and a generalized convolution, \textit{Integr. Trans. Spec. Funct.} \textbf{34} (2023), no. 9, 690–702. \url{https://doi.org/10.1080/10652469.2023.2182776}
\bibitem{Stein1971Weiss}E. M. Stein and G. Weiss, \textit{Introduction to Fourier analysis on Euclidean spaces}, Princeton Mathematical Series (PMS-32), Vol.  32, Princeton University Press, NJ, 1972.
\bibitem{Tsitsiklis1981Levy}J. N. Tsitsiklis and B. C. Levy: Integral equations and resolvents of Toeplitz plus Hankel kernels. Series/ Report. No: LIDS-P 1170. Laboratory for Information and Decision Systems, Massachusetts Institute of Technology, 1981. \url{http://hdl.handle.net/1721.1/1005}
\bibitem{Appell2000KalZabre}J. M. Appell, A. S. Kalitvin, and P. P. Zabrejko, \textit{Partial integral operators and integro-differential equations: Pure and applied mathematics}, 1st ed., Marcel Dekker Inc, New York, 2000. \url{https://doi.org/10.1201/9781482270402}
\bibitem{Agranovich1963Marchenko}Z. S. Agranovich and V. A. Marchenko, \textit{The inverse problem of scattering theory}, Gordon and Breach, New York, 1963.
\bibitem{Chadan1977Sabatier}K. Chadan and P. C. Sabatier, \textit{Inverse problems in quantum scattering theory}, Springer-Verlag, New York, 1977.
\bibitem{kagiwada1974integral}H. H. Kagiwada and R. Kalaba, \textit{Integral equations via imbedding methods}, Applied Mathematics and Computation, No. 6, Addison-Wesley Publishing Company, Reading, 1974.
\bibitem{Gelfand1951Levitan}I. M. Gel'fand and B. M. Levitan, On the determination of a differential equation from its spectral function. (in Russian), \textit{Izv. Akad. Nauk SSSR. Ser. Mat.} \textbf{15} (1951), no. 4, 309–360.
\bibitem{tuan2018hoanghong}T. Tuan, P. V. Hoang, and N. T. Hong, Integral equation of Toeplitz plus Hankel's type and parabolic equation related to the Kontorovich-Lebedev-Fourier generalized convolutions, \textit{Math. Methods Appl. Sci.} \textbf{41} (2018), no. 17, 8171–8181. \url{https://doi.org/10.1002/mma.5279}
\bibitem{NaimarkMA1972} Naimark MA. Normed Algebras, Third edition; Wolters–Noordhoff Series of Monographs and Textbooks on Pure and Applied Mathematics Edition, p. 598. Wolters-Noordhoff Publishing, Groningen (1972). Translated from the second Russian edition by Leo F. Boron.
\bibitem{castro2019}L. P. Castro, R. C. Guerra, and N. M. Tuan, Convolution theorems related with the solvability of Wiener-Hopf plus Hankel integral equations and Shannon's sampling formula, Math. Slovaca. 69 (2019), no. 5, 1149–1164. \url{https://doi.org/10.1515/ms-2017-0297}
\bibitem{castro2020new}L. P. Castro, R. C. Guerra, and N. M. Tuan, New convolutions and their applicability to integral equations of Wiener-Hopf plus Hankel type, \textit{Math. Methods Appl. Sci.} \textbf{43} (2020), no. 7, 4835–4846. \url{https://doi.org/10.1002/mma.6236}
\bibitem{Barbashin1957}E. A. Barbashin, Conditions for invariance of stability of solutions of integro-differential equations (in Russian), \textit{Izv. Vyssh. Uchebn. Zaved. Mat.} \textbf{1} (1957), 25–34.
\bibitem{KilbasSrivastavaTrujillo2006}A. A. Kilbas, H. M. Srivastava, and J. J. Trujillo, \textit{Theory and applications of fractional differential equations}, North-Holland Mathematics Studies, Vol.  \textbf{204}, Elsevier, Amsterdam, 2006.
\end{thebibliography}

\end{document}